%Format: plain

%% arXiv:0812.3258

\input amstex
\documentstyle{amsppt}

\input label.def
\input degt.def
%\input debug.def
%\PrintLabels

%\copycounter\thm\subsection
%\copycounter\equation\subsection

\input epsf
\def\picture#1{\epsffile{#1-bb.eps}}
\def\cpic#1{$\vcenter{\hbox{\picture{#1}}}$}

\def\ie{\emph{i.e.}}
\def\eg{\emph{e.g.}}
\def\cf{\emph{cf}}
\def\via{\emph{via}}
\def\etc{\emph{etc}}

%% From the previous paper %%

%%%%%%%%%%%%%%%%%%%%%%%%%%%%%

{\catcode`\@11
\gdef\proclaimfont@{\sl}}

\newtheorem{Remark}\Remark\thm\endAmSdef
\conjecture\thm\endproclaim
\newhead\subsubsection\subsubsection\endsubhead
\def\paragraph{\subsubsection{}}

\def\dash{\item"\hfill--\hfill"}
\def\Dashes{\widestnumber\item{--}\roster}
\def\endDashes{\endroster}

\loadbold
\def\bA{\bold A}
\def\bD{\bold D}
\def\bE{\bold E}
\def\bJ{\bold J}

\def\tA#1{\smash{\tilde\bA#1}}
\def\tD#1{\smash{\tilde\bD#1}}
\def\tE#1{\smash{\tilde\bE#1}}

\let\Ge\epsilon

\let\splus\oplus

\let\onto\twoheadrightarrow

\def\F{\Bbb F}

\def\CG#1{\Z_{#1}}
\def\BG#1{\Bbb B_{#1}}
\def\DG#1{\Bbb D_{#1}}

\def\TG#1{\Bbb T_{2,#1}}

\def\jB{j_{\B}}
\def\GB{\Sk}

\let\Gs\sigma

\def\SL{\operatorname{\text{\sl SL\/}}}

\def\gcd{\operatorname{g.c.d.}}

\def\CC{\Cal C}

\def\B{\bar B}

\def\Cp#1{\Bbb P^{#1}}
\def\Rp#1{\Bbb P_{\R}^{#1}}
\def\term#1-{$\DG{#1}$-}

\let\Ga\alpha
\let\Gb\beta

\let\Gs\sigma
\let\Gr\rho
\def\1{^{-1}}

\def\ls|#1|{\mathopen|#1\mathclose|}
\let\<\langle
\let\>\rangle

\def\ord{\operatorname{ord}}
\def\Sk#1{\operatorname{Sk}#1}

\let\MB=M

\def\0{{\setbox0\hbox{0}\hbox to\wd0{\hss\rm--\hss}}}

\def\tabstrut{\vrule height9.5pt depth2.5pt}
\def\exstrut{\omit\vrule height2pt\hss\vrule}

\def\GAP{{\tt GAP}}
\def\fref(#1){\rom{(#1)}}
\def\fragment(#1)#2{\ref{fig.e8}\fref(#1)--$#2$}
\def\bifrag(#1)#2{\fragment(#1){#2,\bar#2\!\!\!}}
\def\rfragment(#1)#2{\ref{fig.e8-r}\fref(#1)--$#2$}
\def\rbifrag(#1)#2{\rfragment(#1){#2,\bar#2\!\!\!}}

\def\inserthyphen{\ifcat\next a-\fi\ignorespaces}
\let\BLACK\bullet
\let\WHITE\circ
\def\CROSS{\vcenter{\hbox{$\scriptstyle\mathord\times$}}}
\let\STAR*
\def\pblack-{$\BLACK$\futurelet\next\inserthyphen}
\def\pwhite-{$\WHITE$\futurelet\next\inserthyphen}
\def\pcross-{$\CROSS$\futurelet\next\inserthyphen}
\def\pstar-{$\STAR$\futurelet\next\inserthyphen}
\def\black{\protect\pblack}
\def\white{\protect\pwhite}
\def\cross{\protect\pcross}
\def\star{\protect\pstar}
\def\NO#1{\mathord\#_{#1}}
\def\nblack{\NO\BLACK}
\def\nwhite{\NO\WHITE}

\def\nstar{\NO\STAR}

\let\bb\delta
\let\bc\beta

\def\cstar{$(*)$}

\def\beginGAP{\bgroup
 \catcode`\^=12\catcode`\#=12\catcode`\_=12
 \obeylines\obeyspaces\eightpoint\tt}
\let\endGAP\egroup

\topmatter

\author
Alex Degtyarev
\endauthor

\title
Plane sextics with a type $\bE_8$ singular point
\endtitle

\address
Department of Mathematics,
Bilkent University,
06800 Ankara, Turkey
\endaddress

\email
degt\@fen.bilkent.edu.tr
\endemail

\abstract
We construct explicit geometric models for and
compute the fundamental groups of all plane sextics with simple
singularities only
and with at least one type~$\bE_8$ singular point.
In particular, we discover four new sextics with nonabelian
fundamental groups; two of them are irreducible. The groups of the
two
irreducible sextics found are finite. The principal tool used is
the reduction to trigonal curves and
Grothendieck's \emph{dessins d'enfants}.
\endabstract

\keywords
Plane sextic, singular curve,
fundamental group, trigonal curve, dessin d'enfant
\endkeywords

\subjclassyear{2000}
\subjclass
Primary: 14H45; % curves/Special curves and curves of low genus
Secondary: 14H30, % curves/Coverings, fundamental group
14H50 % curves/Plane and space curves
\endsubjclass

\endtopmatter

\document

\section{Introduction}

\subsection{Principal results}
This paper is a continuation of my paper~\cite{dessin.e7},
where we started
the study of the equisingular deformation families and the
fundamental groups of plane \emph{sextics} (\ie,
%projective plane
curves $B\subset\Cp2$
of degree six) with a
distinguished triple
singular point, using the representation of such sextics
\via\ trigonal curves in Hirzebruch surfaces and Grothendieck's
\emph{dessins d'enfants}. (All varieties
%considered
in the paper
are over~$\C$
and are considered in their Hausdorff topology.)
Recall that, in spite of the fact that
the deformation classification of sextics can be reduced to a
relatively simple, although tedious, arithmetical problem,
see~\cite{JAG}, the geometry of the pairs $(\Cp2,B)$
remains a \emph{terra incognita}, as the construction relies upon
the global Torelli theorem for $K3$-surfaces and is quite
implicit. On the contrary, the approach suggested
in~\cite{dessin.e7}, although not resulting in a defining equation
for~$B$, gives one a fairly good understanding of the topology
of $(\Cp2,B)$; in particular, it is sufficient for the computation
of the fundamental group $\pi_1(\Cp2\sminus B)$. A few other
applications of this approach
%are discussed
and more motivation can be found
in~\cite{dessin.e7};
for a brief overview of the latest achievements on the subject,
%we refer to
see
C.~Eyral, M.~Oka~\cite{EyralOka.new}.

In the present paper, we deal with the case when the
distinguished triple point in question is of
type~$\bE_8$. (The case of a type~$\bE_7$ singular point was
considered in~\cite{dessin.e7}, and the case of~$\bE_6$ is the
subject of a forthcoming paper.) As in~\cite{dessin.e7}, a simple
trick with the skeletons reduces most sextics $B\subset\Cp2$
to certain
trigonal curves~$\B'$
in~$\Sigma_2$ (instead of the original surface~$\Sigma_3$);
this simplifies dramatically the classification of the sextics
and the computation of their fundamental
groups.
%The geometric relation between~$B$ and~$\B'$ remains a
%mystery.
It is still unclear if there is a simple
geometric relation between~$B$ and~$\B'$.
%Unfortunately, their seems to be no similar construction
%for the type~$\bE_6$ singular point.

Throughout the paper, we consider a plane sextic $B\subset\Cp2$
satisfying the following conditions:
%\widestnumber\item{$(*)$}
\roster
\item"\cstar"
$B$ has simple singularities only, and
\item""
$B$ has a distinguished singular point~$P$ of type~$\bE_8$.
\endroster
(It is worth mentioning that all sextics
with a non-simple singular point,
as well as their fundamental groups, are well known.
If such a sextic~$B$ also
has a type~$\bE_8$ singular point, then the equisingular
deformation class of~$B$ is determined by its set of singularities,
which is either $\bE_8\splus\bE_{12}$ or $\bE_8\splus\bJ_{2,i}$,
$i=0,1$, in Arnol$'$d's notation,
and
the fundamental group $\pi_1(\Cp2\sminus B)$ is abelian.)

Recall that the total Milnor number $\mu(B)$ of a plane sextic
$B\subset\Cp2$
with simple singularities only is subject to the inequality
$\mu(B)\le19$. (The double covering~$X$ of
the plane~$\Cp2$ ramified at~$B$
is a $K3$-surface, and since the exceptional divisors arising from
the resolution of the singular points are all linearly independent
and span a negative definite lattice,
one has $\mu(B)\le h^{1,1}(X)-1=19$.) A
sextic~$B$ is called \emph{maximal} (sometimes \emph{maximizing})
if $\mu(B)=19$.
Recall also that maximal sextics are \emph{rigid}, \ie, two such
sextics belong to a connected equisingular deformation family if
and only if they are related by a projective transformation.
The principal results of the present paper are
Theorems~\ref{th.e8} and~\ref{th.e8-r} (an explicit
classification of the maximal sextics with a type~$\bE_8$ singular
point) and Theorems \ref{th.main} and~\ref{th.pert} (the
computation of the fundamental groups).

\theorem\label{th.e8}
Up to projective transformation \rom(equivalently, up to
equisingular deformation\rom),
there are $39$ maximal irreducible
sextics~$B$ satisfying~\cstar\rom;
they realize $26$ sets of
singularities \rom(see Table~\ref{tab.e8} on
Page~\rom{\pageref{tab.e8}}\rom).
\endtheorem

\theorem\label{th.e8-r}
Up to projective transformation \rom(equivalently, up to
equisingular deformation\rom),
there are $18$ maximal reducible
sextics~$B$ satisfying~\cstar\rom;
they realize $17$ sets of
singularities \rom(see Table~\ref{tab.e8-r} on
Page~\rom{\pageref{tab.e8-r}}\rom).
\endtheorem

Theorems~\ref{th.e8} and~\ref{th.e8-r}
are proved in~\S\ref{S.classification}
(see~\ref{proof.e8} and~\ref{proof.e8-r}, respectively), where
all
sextics are constructed explicitly using trigonal curves.
(This construction is
further
used in~\S\ref{S.groups} in the computation
of the fundamental groups of the sextics.)
Alternatively, the statements
%of the theorems
could as well be derived
from combining
the results of J.-G.~Yang~\cite{Yang} (the existence) and
I.~Shimada~\cite{Shimada} (the enumeration of the maximal
sets of
singularities realized by more than one deformation family).

\theorem\label{th.main}
With two exceptions, the fundamental group $\pi_1(\Cp2\sminus B)$
of a
%maximal
plane sextic $B\subset\Cp2$
satisfying~\cstar\ is abelian.
The exceptions are\rom:
\roster
\item\local{G.irr}
the
%\rom(only\rom)
sextic with the set
of singularities $\bE_8\splus\bA_4\splus\bA_3\splus2\bA_2$\rom;
%its group is given by
the group is
$$
G_6:=\bigl\<\Ga_1,\Ga_2\bigm|
 (\Ga_1\Ga_2\1)^5\Ga_2^6=1,\
 [\Ga_1,\Ga_2^3]=1,\
 (\Ga_1\Ga_2)^2\Ga_1=(\Ga_2\Ga_1)^2\Ga_2\bigr\>;
$$
\item\local{G.red}
the
%\rom(only\rom)
sextic with the set
of singularities $\bE_8\splus\bD_6\splus\bA_3\splus\bA_2$\rom;
%its group is given by
the group is
$$
G_\infty:=\bigl\<\Ga_1,\Ga_2\bigm|
 (\Ga_1\Ga_2\1)^5\Ga_2^6=1,\
 [\Ga_1,\Ga_2^3]=1,\
 (\Ga_1\Ga_2)^2=(\Ga_2\Ga_1)^2\bigr\>.
$$
\endroster
Each group~$G_i$ above, $i=6$ or~$\infty$, can be represented as a
semi-direct
product of its abelianization, which is a cyclic group
of order~$i$, and its commutant
$[G_i,G_i]\cong\SL(2,\Bbb F_5)$, which is the only perfect group of
order $120$.
\endtheorem

Theorem~\ref{th.main} is proved in~\ref{proof.main}. It is worth
mentioning that the group~$G_6$ in \iref{th.main}{G.irr} is
finite, substantiating my conjecture that the group of \emph{any}
irreducible sextic with simple singularities that is not of torus
type is finite. (Recall that a plane sextic~$B$ is said to be of
\emph{torus type} if its equation can be represented in the form
$f_2^3+f_3^2=0$ for some homogeneous polynomials~$f_2$, $f_3$ of
degree~$2$ and~$3$, respectively. The groups of sextics of torus
type are all infinite: they factor ro the reduced braid group
$\BG3/(\Gs_1\Gs_2)^3\cong\CG2*\CG3$.)

\Remark\label{rem.groups}
More precisely, one has
$$
G_\infty=\Z\times\SL(2,\F_5)
\quad\text{and}\quad
G_6=\CG{12}\odot\SL(2,\F_5),
$$
where the latter central product is the quotient of
$\CG{12}\times\SL(2,\F_5)$
by the diagonal
$\CG2\subset\CG2\times\operatorname{Center}\SL(2,\F_5)$;
for details, see~\ref{s.120}
and~\ref{s.720}, respectively.
\endRemark

Recall that, due to O.~Zariski~\cite{Zariski.group}, any
perturbation $B\to B'$ of reduced plane curves induces an
epimorphism $\pi_1(\Cp2\sminus B)\onto\pi_1(\Cp2\sminus B')$
of their fundamental groups. In particular, if
$\pi_1(\Cp2\sminus B)$ is abelian, so is $\pi_1(\Cp2\sminus B')$.
Next theorem describes the few perturbations of plane sextics
as in Theorem~\ref{th.main} that have nonabelian
fundamental groups.

\theorem\label{th.pert}
With two exceptions, the fundamental group of a sextic~$B'$
that is a proper
perturbation of a plane sextic~$B$
%as in
%Theorem~\ref{th.main}
satisfying~\cstar\
is abelian. The exceptions are\rom:
\roster
\item
the perturbation
$\bE_8\splus\bA_4\splus\bA_3\splus2\bA_2\to
2\bA_4\splus2\bA_3\splus2\bA_2$, and
\item
the perturbation
$\bE_8\splus\bD_6\splus\bA_3\splus\bA_2\to
\bD_6\splus\bD_5\splus\bA_3\splus2\bA_2$.
\endroster
For both curves, the perturbation epimorphism
$\pi_1(\Cp2\sminus B)\onto\pi_1(\Cp2\sminus B')$
is an isomorphism\rom;
%hence, the
in particular, the fundamental
groups of the curves are $G_6$ and~$G_\infty$,
respectively, see Theorem~\ref{th.main}.
\endtheorem

This theorem is proved in~\ref{proof.pert}. It covers over two
hundred new (compared to~\cite{dessin.e7})
sets of simple singularities realized by sextics with
abelian fundamental groups.
(Recall that, according to~\cite{degt.8a2},
any induced subgraph of the
combined Dynkin graph of a plane sextic~$B$ with simple
singularities only
can be realized by a perturbation of~$B$; in other words, the
singular points of~$B$ can be perturbed independently.)
The total number of such sets of
singularities currently known is over $1400$.

\subsection{Classical Zariski pairs}\label{s.Zariski}
Among the new sets of singularities
realized by sextics with abelian fundamental groups
is
$\bE_6\splus\bA_8\splus\bA_2\splus2\bA_1$ (a perturbation of
$\bE_8\splus\bA_8\splus\bA_2\splus\bA_1$,
Nos\.~$9$ and~$17$ in Table~\ref{tab.e8}). The corresponding
sextic is included into a so called \emph{classical Zariski pair},
\ie, a pair of irreducible sextics that share the same set of
singularities but
%differ by their
have different
Alexander polynomials (see,
\eg,~\cite{JAG} for details). In each pair, one of the curves is
\emph{abundant}, or of \emph{torus type}, and its Alexander
polynomial is $t^2-t+1$; the other curve is \emph{non-abundant},
or not of torus type, and its Alexander polynomial is~$1$.
(The term `abundant' is due to the fact that the Alexander
polynomial of the curve is larger than the minimal polynomial
imposed by its singularities.)
Conjecturally, in each pair the fundamental group of the abundant
curve is the reduced braid group $\BG3/(\Gs_1\Gs_2)^3\cong\CG2*\CG3$,
whereas
the group of the non-abundant curve is abelian
(hence, equal to~$\CG6$).
At present, the
conjecture is known for all sets of singularities except
$$
\bA_{17}\splus\bA_1,\quad
\bA_{14}\splus\bA_2\splus2\bA_1,\quad %(torus type is unknown either),
2\bA_8\splus2\bA_1,\quad
2\bA_8\splus\bA_1.
$$
(The group of the abundant curve is known for all sets of
singularities except $\bA_{14}\splus\bA_2\splus2\bA_1$.) In
E.~Artal \emph{et al\.}~\cite{Artal.Trends}, it is stated that the
fundamental group of a sextic
with a single type~$\bA_{19}$ singular point is abelian;
by perturbation, this assertion implies that the sets of
singularities $\bA_{17}\splus\bA_1$
and $2\bA_8\splus\bA_1$
%are
%also
can also be
realized by irreducible sextics
with abelian groups, leaving the conjecture unsettled for two sets
of singularities only.

\subsection{Contents of the paper}\label{s.contents}
The paper depends on a preliminary computation found
in~\cite{dessin.e7}; it is based on the theory of trigonal curves,
Grothendieck's \emph{dessins d'enfants}, braid monodromy, and
Zariski--van Kampen's method. We refer to~\cite{dessin.e7} for a
brief exposition of these subjects.

In~\S\ref{S.classification}, we prove Theorems~\ref{th.e8}
and~\ref{th.e8-r} by providing an explicit geometric construction,
in terms of the skeletons of the trigonal models,
for all maximal sextics satisfying~\cstar. This construction is
used in~\S\ref{S.groups} to compute the fundamental groups of the
maximal sextics. In~\S\ref{S.perturbations}, we analyze the
perturbations of
%the
maximal sextics
with a type~$\bE_8$ singular point
and prove
Theorems~\ref{th.main} and~\ref{th.pert}. An important technical
result here is Proposition~\ref{pert.e8}, describing the
local fundamental groups of all perturbations of a type~$\bE_8$
%singular point.
singularity.

\section{The classification}\label{S.classification}

In this section, we prove Theorems~\ref{th.e8} and~\ref{th.e8-r}.

\subsection{The trigonal model}\label{s.model}
The following statement, proved in~\cite{dessin.e7}, reduces the
study of sextics with a type~$\bE_8$ singular point to the study
of trigonal curves in the Hirzebruch surface~$\Sigma_3$. (Recall
that $\Sigma_3$ is a geometrically ruled rational surface with an
exceptional section~$E$ if self-intersection~$-3$; we refer
to~\cite{dessin.e7} and~\cite{degt.kplets} for the terminology and
notation related to trigonal curves.)

\proposition\label{1-1.e8}
There is a natural
bijection~$\phi$, invariant under equisingular deformations,
between the following two sets\rom:
\roster
\item\local{e8.sextic}
plane sextics~$B$ with a distinguished
type~$\bE_8$ singular point~$P$, and
\item\local{e8.trigonal}
trigonal curves $\B\subset\Sigma_3$ with a
distinguished type~$\tA{_1^*}$ singular fiber~$F$.
\endroster
A sextic~$B$ is irreducible if and only if so is $\B=\phi(B)$ and,
with one exception, $B$ is maximal if and only if $\B$ is
maximal and has no unstable fibers other than~$F$.
The exception is the reducible sextic~$B$ with the set of
singularities $\bE_8\splus\bE_7\splus\bD_4$\rom; in this case,
$\phi(B)$ is isotrivial.
\endproposition

The trigonal curve~$\B$ corresponding to a sextic~$B$ (more
precisely, pair $(B,P)$) under~$\phi$ above
is called the \emph{trigonal model}
of~$B$.

From now on, we disregard the exceptional case
$\bE_8\splus\bE_7\splus\bD_4$ and assume that $\B$ is not
isotrivial, hence maximal. Let $\jB\:\Cp1\to\Cp1$ be the
$j$-invariant of~$\B$, and let $\GB=\jB\1([0,1])\subset S^2$ be
its skeleton, see~\cite{dessin.e7} or~\cite{degt.kplets}. Recall
that $\GB$ is a connected planar map with at most $3$-valent
\black-vertices (the pull-backs of~$0$)
and monovalent \white-vertices (some of the
pull-backs of~$1$) connected to
\black-vertices. (In fact, each edge connecting two
\black-vertices contains a bivalent \white-vertex in the middle,
making~$\GB$ a bipartite graph,
but these bivalent \white-vertices are ignored.)

Denote by~$u$ the monovalent
\white-vertex of~$\GB$ corresponding to
the distinguished type~$\tA{_1^*}$ fiber~$F$ given by
Proposition~\ref{1-1.e8}, and let~$v$ be the \black-vertex adjacent
to~$u$. The induced subgraph of~$\GB$
spanned by~$u$ and~$v$ (\ie, $u$, $v$, and the edge $[u,v]$)
is called the
\emph{insertion}; in the drawings below it is shown in grey.

\lemma\label{type.E}
The mono-
\rom(respectively, bi-\rom) valent \black-vertices of~$\GB$
correspond to the type~$\tE{_6}$
\rom(respectively,~$\tE{_8}$\rom)
singular fibers of~$\B$, and the monovalent \white-vertices
of~$\Sk$ other than~$u$ correspond to the type~$\tE{_7}$ singular
fibers of~$\B$. Furthermore, one has
$$
3\nblack(1)+4\nblack(2)+\nblack(3)+3\nwhite(1)=8-2d,
\eqtag\label{vertex.count}
$$
where $\nstar(i)$ is the number of \star-vertices of valency~$i$,
$\STAR=\BLACK$ or $\WHITE$, and $d$ is the number of the
$\tD{}$-type fibers of~$\B$.
\endlemma

\proof
The first statement follows from the fact that
all singular fibers of~$\B$ other than~$F$ are stable and the
relation between the vertices of~$\Sk$ and the fibers of~$\B$,
see~\cite{dessin.e7} or~\cite{degt.kplets}. Then the number~$t$ of
the triple points of~$\B$ is given by
$$
t=d+\nblack(1)+\nblack(2)+\nwhite(1)-1,
$$
and \eqref{vertex.count} follows from the vertex count given by
Corollary~2.5.5 in~\cite{dessin.e7}.
\endproof

\Remark\label{rem.AD}
Recall that, in addition to the $\tE{}$-type singular
fibers described in
Lemma~\ref{type.E}, each $m$-gonal region of~$\GB$ contains a
unique type~$\tA{_{m-1}}$ ($\tA{_0^*}$ if $m=0$)
or type~$\tD{_{m+4}}$ fiber of~$\B$. The total number~$d$ of the
$\tD{}$-type fibers is given by~\eqref{vertex.count}, and the $d$
regions containing such fibers can be chosen arbitrary, resulting
in general
in distinct deformation classes of curves (and even in distinct
sets of singularities).
\endRemark

\subsection{The reduction}
We still assume that the trigonal model $\B\subset\Sigma_3$
is non-isotrivial and maximal and keep
the notation~$u$ and~$v$ for the \white-- and, respectively,
\black-vertices spanning the insertion in~$\GB$.

\lemma\label{val.1}
If $v$ is a monovalent vertex, then $\GB$ is as shown in
Figure~\ref{fig.e8}\rom{(k)} below and the set of singularities of~$B$ is
$\bE_8\splus\bE_6\splus\bD_5$.
\endlemma

\proof
Under the assumptions,
the insertion is an edge bounded by two monovalent vertices. On
the other hand, it is a subgraph of a connected graph~$\GB$.
Hence, $\GB$ is exhausted by the insertion, \ie, it is the graph
shown in Figure~\ref{fig.e8}\rom{(k)}.
Then, \eqref{vertex.count}
implies that $d=1$; hence, the set of singularities of~$B$ is
$\bE_8\splus\bE_6\splus\bD_5$.
\endproof

\lemma\label{val.2}
If $v$ is a bivalent vertex, then $\GB$ is the graph shown in
Figure~\ref{fig.e8}\rom{(j)} below and the set of
singularities of~$B$ is
$2\bE_8\splus\bA_3$.
\endlemma

\proof
If $v$ is bivalent (corresponding to a type~$\tE{_8}$ singular
fiber of~$\B$, \ie, to
%another
a
type~$\bE_8$ singular point
of~$B$
other than~$P$), then the vertex count~\eqref{vertex.count}
implies that $\GB$ has exactly one more
\black-vertex, which is trivalent, and one can easily see that the
only skeleton with these properties is the one shown in
Figure~\ref{fig.e8}\fref(j). Besides, \eqref{vertex.count} also implies
$d=0$. Hence, the set of singularities of~$B$ is
$2\bE_8\splus\bA_3$.
\endproof

Finally, consider the `general' case, when $v$ is trivalent. In
this case, removing the insertion and
patching the two other edges incident to~$v$ to a single edge, one
obtains a new connected graph~$\Sk'$.

\lemma\label{Sk'}
With one exception \rom(see below\rom),
the graph~$\Sk'$
constructed
%described
above is the skeleton of a certain
stable maximal trigonal curve $\B'\subset\Sigma_2$.
Conversely,
attaching
an insertion at the middle of any edge of
the skeleton~$\Sk'$ of a
stable maximal trigonal curve $\B'\subset\Sigma_2$, one obtains a
skeleton~$\GB$ defining a maximal trigonal curve~$\B$ as in
Proposition~\ref{1-1.e8}.

The exception
is the skeleton $\GB$ shown in Figure~\ref{fig.e8-r}\rom{(g)},
when $\Sk'$ has no \black-vertices\rom; the set of singularities
of the corresponding curve~$B$ is $\bE_8\splus\bD_6\splus\bD_5$.
\endlemma

\proof
First, notice that $v$ cannot be adjacent to three monovalent
\white-vertices, as that would contradict to~\eqref{vertex.count}.
Then, it is easy to see that,
provided that $\Sk'$ has at least one \black-vertex
(as any skeleton does),
%(which is a
%common property of all skeletons),
it is indeed a valid skeleton,
and \eqref{vertex.count}
transforms to the following vertex count for~$\Sk'$:
$$
3\nblack(1)+4\nblack(2)+\nblack(3)+3\nwhite(1)=4-2d.
$$
Using Corollary~2.5.5
in~\cite{dessin.e7}, one concludes that $\Sk'$ is the skeleton of a
stable maximal trigonal curve $\B'\subset\Sigma_2$.
The converse statement is obvious.

In the exceptional case, when $\Sk'$ has no \black-vertices, since
$\Sk'$ is still connected, it must consist of a single circle;
then $\GB$ is the graph shown in Figure~\ref{fig.e8-r}\rom{(g)},
and the vertex count~\eqref{vertex.count} implies that $d=2$, \ie,
each of the two regions of~$\Sk$ contains a $\tD{}$-type fiber
of~$\B$ and the set of
singularities of~$B$ is $\bE_8\splus\bD_6\splus\bD_5$.
\endproof

\subsection{Reducible vs\. irreducible curves}
Recall that a \emph{marking} at a trivalent \black-vertex~$w$
of a skeleton~$\Sk$
is a counterclockwise ordering $\{e_1,e_2,e_3\}$
of the three edges
attached to~$w$.
A marking is uniquely defined by assigning index~$1$ to one of
the three edges.
Given a marking,
the indices of the edges are considered defined
modulo~$3$, so that $e_4=e_1$, $e_5=e_2$, \etc.

%A marking at~$v$ defines as well
%an ordering $\{e'_1,e'_2,e'_3\}$ of the three
%solid edges attached to~$v$: we let~$e'_i$
%to be the solid edge opposite to~$e_i$.

\definition\label{def.splitting}
A \emph{marking} of a
skeleton~$\Sk$ is a collection of markings at all trivalent
\black-vertices of~$\Sk$. Given a marking, one can assign a type $[i,j]$,
$i,j\in\Z_3$, to each edge~$e$ of~$\Sk$ connecting two trivalent
\black-vertices, according to the indices of the two ends of~$e$. A
marking of a skeleton without mono- or bivalent \black-vertices
is called \emph{splitting} if it satisfies the following
two conditions:
\roster
\item\local{3.types}
the types of all edges are $[1,1]$, $[2,3]$, or~$[3,2]$;
\item\local{[bw]}
an edge connecting a \black-vertex~$w$ and a monovalent
\white-vertex
%has index~$1$
is~$e_1$
at~$w$.
\endroster
\enddefinition

According to~\cite{degt.kplets}, the splitting markings of the
skeleton of a maximal trigonal curve $\B\subset\Sigma_k$
are in a one-to-one
correspondence with the \emph{linear}
components of~$\B$ (\ie, components of~$\B$ that are sections
of~$\Sigma_k$).

\lemma\label{prop.reducible}
Let
%$\B$, $\GB$ and $\B'$, $\Sk'$ be
$\Sk$, $\Sk'$ be a pair of skeletons
as in Lemma~\ref{Sk'}, so
that $\Sk'$ is obtained from~$\GB$ by removing the insertion.
%Let~$\Gamma_0$ be a skeleton, and let $\B$ be the trigonal curve
%whose skeleton is obtained from~$\Gamma_0$ by adding an insertion.
Then
the curve~$\B$ defined by~$\GB$
is reducible if and only if $\Sk'$ admits a
splitting marking such that the insertion is attached to an edge of
type $[2,3]$ and is oriented as shown in Figure~\ref{fig.ins}{\rm(b)}.
\endlemma

\midinsert
\centerline{\vbox{\halign{\hss#\hss&&\kern6em\hss#\hss\cr
\cpic{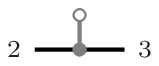}&\cpic{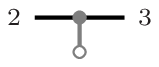}\cr
\noalign{\medskip}
(a)&(b)\cr
\crcr}}}
\figure
An insertion resulting in (a) irreducible and (b) reducible
trigonal curve~$\B$
\endfigure\label{fig.ins}
\endinsert

\proof
%Using condition~\iref{def.splitting}{[bw]}, one can see
%that
In view of condition~\iref{def.splitting}{[bw]},
a splitting marking of~$\GB$ restricts to a splitting
marking of~$\Sk'$. Conversely,
due to~\iref{def.splitting}{[bw]} again, a
splitting marking of~$\Sk'$ extends to~$\GB$ if and only if the
insertion is as shown in Figure~\ref{fig.ins}\fref(b).
\endproof

\corollary\label{2.markings}
In the notation of Lemma~\ref{prop.reducible}, if $\Sk'$
admits more than one splitting marking \rom(equivalently, if
the corresponding trigonal curve~$\B'$
splits into three distinct linear components\rom),
then $\B$ is reducible.
\endcorollary

\proof
A splitting marking is uniquely determined by its restriction to a
single vertex. Hence, the restrictions of the three splitting
markings of~$\Sk'$ to any given \emph{oriented} edge~$e$ are pairwise
distinct and, since there are only three types, see
condition~\iref{def.splitting}{3.types},
$e$
%is of type $[2,3]$
can be made of any given type
under one of the markings.
\endproof

\midinsert
\centerline{\vbox{\halign{\hss#\hss&&\qquad\hss#\hss\cr
\cpic{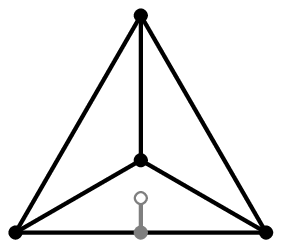}&\cpic{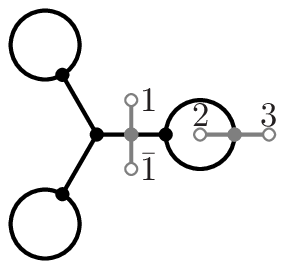}&\cpic{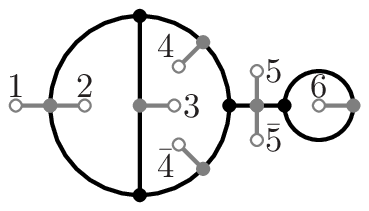}\cr
\noalign{\medskip}
(a)&(b)&(c)\cr
\crcr}}}
\bigskip
\centerline{\vbox{\halign{\hss#\hss&&\qquad\quad\hss#\hss\cr
\cpic{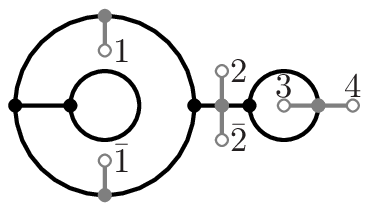}&\cpic{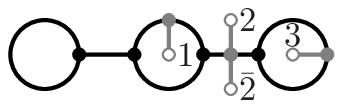}\cr
\noalign{\medskip}
(d)&(e)\cr
\crcr}}}
\bigskip
\centerline{\vbox{\halign{\hss#\hss&&\kern1.3em \hss#\hss\cr
\cpic{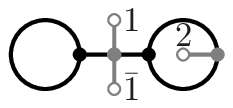}&\cpic{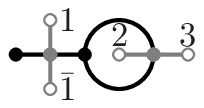}&\cpic{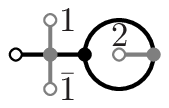}&\cpic{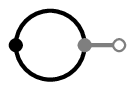}&\cpic{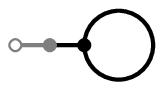}&\cpic{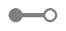}\cr
\noalign{\medskip}
(f)&(g)&(h)&(i)&(j)&(k)\cr
\crcr}}}
\figure\label{fig.e8}
Skeletons of irreducible curves~$\B$
\endfigure
\endinsert

\subsection{Proof of Theorem~\ref{th.e8}}\label{proof.e8}
Due to Proposition~\ref{1-1.e8}, we need to enumerate all
irreducible
maximal trigonal curves $\B\subset\Sigma_3$ with a unique unstable
fiber, which is of type~$\tA{_1^*}$. The study of such curves
reduces, in its turn, to the study of the skeletons containing an
insertion and satisfying the vertex count~\eqref{vertex.count}.
With two exceptions, described in Lemmas~\ref{val.1}
and~\ref{val.2} (see Figure~\ref{fig.e8}\fref(j) and~\fref(k);
these two
skeletons obviously define irreducible curves), each skeleton
$\Sk$ as above is obtained by attaching an insertion to the
skeleton~$\Sk'$ of a stable maximal trigonal curve
$\B'\subset\Sigma_2$, see Lemma~\ref{Sk'}. (It is easy to see
that the exceptional
skeleton in Lemma~\ref{Sk'} admits a splitting marking and thus
defines a reducible curve.)

\midinsert
\table\label{tab.e8}
Maximal sets of singularities with a type~$\bE_8$ point
represented by irreducible sextics
\endtable

\def\no{}
\let\-\0
\def\*{\llap{$^*$}}
\let\*\relax
\def\tref#1{\,\text{\ref{#1}}\,}
\def\same#1{\llap{$^{#1\,}$}}
\centerline{\vbox{\offinterlineskip\halign{%
\tabstrut\ \hss#\hss\ \vrule&
\quad$#$\hss\quad\vrule&
\quad\null#\hss\quad\vrule&&\ \ \hss$#$\hss\ \ \vrule\cr
\noalign{\hrule}
\exstrut&&&&&\cr
\#&\text{Set of singularities}&
 \hss Figure&\text{Count}&
 %\text{Real?}&
 \pi_1&
 (l,m,n)
 %S^\perp
 \cr
\exstrut&&&&&\cr
\noalign{\hrule}
\exstrut&&&&&\cr
1&\bE_8\splus\bA_4\splus\bA_3\splus2\bA_2&
 \ref{fig.e8}\fref(a)&(1,0)&\tref{s.720}&
 \*(5,4,3)
 %(15,0,6)
 \cr
2&\bE_8\splus\bA_{11}&
 \bifrag(b){1}&(0,1)&\tref{s.stem}&
 (\-,\-,1)
 %(3,0,2)
 \cr
3&\same1 \bE_8\splus\bA_9\splus\bA_2&
 \fragment(b)2&(1,0)&\tref{s.leaf-in}&
 (3,\-,\-)
 %(5,0,3)
 \cr
4&\same2 \bE_8\splus\bA_{10}\splus\bA_1&
 \fragment(b)3&(1,0)&\tref{s.leaf-out}&
 %(11,0,1)
 \cr
5&\bE_8\splus\bA_7\splus2\bA_2&
 \fragment(c)1&(1,0)&\tref{s.others}&
 (8,3,3)
 %(12,0,3)
 \cr
6&\bE_8\splus\bA_6\splus\bA_3\splus\bA_2&
 \fragment(c)2&(1,0)&\tref{s.others}&
 (4,7,3)
 %(21,0,2)
 \cr
7&\same3 \bE_8\splus\bA_5\splus\bA_4\splus\bA_2&
 \fragment(c)3&(1,0)&\tref{s.others}&
 (5,3,6)
 %(15,0,3)
 \cr
8&\same4 \bE_8\splus\bA_6\splus\bA_4\splus\bA_1&
 \bifrag(c)4&(0,1)&\tref{s.others}&
 (5,7,2)
 %(9,2,4)
 \cr
9&\same5 \bE_8\splus\bA_8\splus\bA_2\splus\bA_1&
 \bifrag(c)5&(0,1)&\tref{s.stem}&
 (\-,\-,1)
 %(9,0,3)
 \cr
10&\bE_8\splus\bA_6\splus2\bA_2\splus\bA_1&
 \fragment(c)6&(1,0)&\tref{s.leaf-in}&
 (3,\-,\-)
 %(21,0,3)
 \cr
11&\bE_8\splus\bA_6\splus\bA_5&
 \bifrag(d)1&(0,1)&\tref{s.others}&
 (7,\-,\-)
 %(11,2,2)
 \cr
12&\bE_8\splus\bA_7\splus\bA_4&
 \bifrag(d)2&(0,1)&\tref{s.stem}&
 (\-,\-,1)
 %(7,2,3)
 \cr
13&\same3 \bE_8\splus\bA_5\splus\bA_4\splus\bA_2&
 \fragment(d)3&(1,0)&\tref{s.leaf-in}&
 (3,\-,\-)
 %(15,0,3)
 \cr
14&\same4 \bE_8\splus\bA_6\splus\bA_4\splus\bA_1&
 \fragment(d)4&(1,0)&\tref{s.leaf-out}&
 %(35,0,1)
 \cr
15&\bE_8\splus\bA_8\splus\bA_3&
 \fragment(e)1&(1,0)&\tref{s.others}&
 (4,9,9)
 %(9,0,2)
 \cr
16&\same2 \bE_8\splus\bA_{10}\splus\bA_1&
 \bifrag(e)2&(0,1)&\tref{s.stem}&
 (\-,\-,1)
 %(4,2,3)
 \cr
17&\same5 \bE_8\splus\bA_8\splus\bA_2\splus\bA_1&
 \fragment(e)3&(1,0)&\tref{s.leaf-in}&
 (3,\-,\-)
 %(9,0,3)
 \cr
18&\bE_8\splus\bD_{11}&
 \fragment(f)1&(1,0)&\tref{s.stem}&
 (\-,\-,1)
 %(2,0,1)
 \cr
19&\bE_8\splus\bD_5\splus\bA_6&
 \bifrag(f)1&(0,1)&\tref{s.stem}&
 (\-,\-,1)
 %(5,2,3)
 \cr
20&\bE_8\splus\bD_9\splus\bA_2&
 \fragment(f)2&(1,0)&\tref{s.leaf-in}&
 (3,\-,\-)
 %(3,0,2)
 \cr
21&\bE_8\splus\bD_7\splus\bA_4&
 \fragment(f)2&(1,0)&\tref{s.others}&
 (\-,5,5)
 %(10,0,1)
 \cr
22&\bE_8\splus\bD_5\splus\bA_4\splus\bA_2&
 \fragment(f)2&(1,0)&\tref{s.leaf-in}&
 (3,\-,\-)
 %(10,0,3)
 \cr
23&\bE_8\splus\bE_6\splus\bA_5&
 \bifrag(g)1&(0,1)&\tref{s.stem}&
 (\-,\-,1)
 %(3,0,3)
 \cr
24&\bE_8\splus\bE_6\splus\bA_3\splus\bA_2&
 \fragment(g)2&(1,0)&\tref{s.leaf-in}&
 (3,\-,\-)
 %(6,0,3)
 \cr
25&\bE_8\splus\bE_6\splus\bA_4\splus\bA_1&
 \fragment(g)3&(1,0)&\tref{s.leaf-out}&
 %(15,0,1)
 \cr
26&\bE_8\splus\bE_7\splus\bA_4&
 \bifrag(h)1&(0,1)&\tref{s.stem}&
 (\-,\-,1)
 %(3,2,2)
 \cr
27&\bE_8\splus\bE_7\splus2\bA_2&
 \fragment(h)2&(1,0)&\tref{s.leaf-in}&
 (3,\-,\-)
 %(3,0,3)
 \cr
28&2\bE_8\splus\bA_2\splus\bA_1&
 \ref{fig.e8}\fref(i)&(1,0)&\tref{s.leaf-out}&
 %(3,0,1)
 \cr
29&2\bE_8\splus\bA_3&
 \ref{fig.e8}\fref(j)&(1,0)&\tref{s.2e8+a3}&
 \*
 %(2,0,1)
 \cr
30&\bE_8\splus\bE_6\splus\bD_5&
 \ref{fig.e8}\fref(k)&(1,0)&\tref{s.e8+e6+d5}&
 \*
 %(6,0,1)
 \cr
\exstrut&&&&&\cr
\noalign{\hrule}
\crcr}}}
\endinsert

The classification of stable maximal
trigonal curves in~$\Sigma_2$ is found in~\cite{symmetric}.
In view of Corollary~\ref{2.markings}, one should take for~$\Sk'$
a skeleton admitting none
(Figure~\ref{fig.e8}\fref(a), \fref(b), \fref(d), \fref(g),
and~\fref(i)) or exactly one
(Figure~\ref{fig.e8}\fref(c), \fref(e), \fref(f),
and~\fref(h)) splitting marking. In the former case, the insertion can
be attached arbitrarily at the middle of any edge; in the latter
case, Lemma~\ref{prop.reducible} implies that
the edge~$e$ that the insertion is attached to is either of
type $[1,1]$ (with respect to the only splitting marking) or of
type $[2,3]$ and oriented with respect to the insertion as shown
in Figure~\ref{fig.ins}\fref(a).
In Figure~\ref{fig.e8}, we list and
number the resulting possibilities for~$\Sk'$
(shown in black) and the
attaching of the insertion (shown in grey), up to symmetries
of~$\Sk'$ (\ie, orientation preserving auto-diffeomorphisms of the
sphere preserving~$\Sk'$).
A pair of indices
$n$, $\bar n$ designates a pair of insertions
that differ by an orientation reversing automorphism of~$\Sk'$;
they result in a pair of complex conjugate sextics.

The set of singularities of~$B$ is almost determined by the
skeleton~$\GB$, see
Lemma \ref{type.E} and Remark~\ref{rem.AD}. (The only
indeterminacy is the one caused by the choice of the region
containing a $\tD{}$-type singular fiber, see
Remark~\ref{rem.AD}.)
The configurations obtained are listed in Table~\ref{tab.e8}. The
table also contains a reference to the fragment in
Figure~\ref{fig.e8} representing the curve and the number of
equisingular deformation classes of curves (in the form
$(n_r,n_c)$, where $n_r$ is the number of real curves and $n_c$ is
the number of pairs of complex conjugate curves, so that the total
number is $n_r+2n_c$). The last two columns are related to the
computation of the fundamental group: the column~`$\pi_1$' refers
to the section where the group is computed, and $(l,m,n)$ are the
values of the parameters used in the computation,
see~\ref{s.others} for details.

Note that some sets of singularities appear from several distinct
skeletons. In Table~\ref{tab.e8}, we prefix the corresponding
lines with an index~$1$, $2$, \etc., equal indices corresponding
to the same set of singularities. (The set
of singularities prefixed with~$1$ is also realized by a reducible
sextic, see Table~\ref{tab.e8-r}.) Summarizing, one obtains $26$
sets of singularities realized by $39$ sextics ($21$ real ones and
$9$ pairs of complex conjugate ones). This proves
Theorem~\ref{th.e8}.
\qed

\Remark
The skeleton~$\Sk'$ in
Figure~\ref{fig.e8}\fref(f) has a symmetry
interchanging fragments~$1$ and~$\bar1$. However, this symmetry
also interchanges the
two monogonal regions of the corresponding skeleton~$\Sk$. Hence,
if one of these regions
contains a type~$\tD{_5}$ fiber of~$\B$
(the set of singularities $\bE_8\splus\bD_5\splus\bA_6$, see
No\.~19 in Table~\ref{tab.e8}), this symmetry does not lift to a
symmetry of $\B\subset\Sigma_3$ and one obtains
two distinct complex conjugate
families.
\endRemark

\midinsert
\centerline{\vbox{\halign{\hss#\hss&&\qquad\qquad\hss#\hss\cr
\cpic{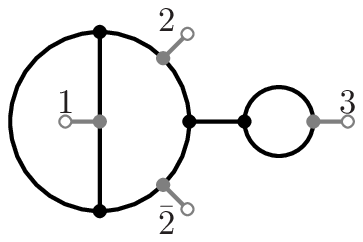}&\cpic{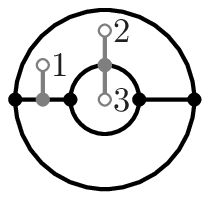}\cr
\noalign{\medskip}
(a)&(b)\cr
\crcr}}}
\bigskip
\centerline{\vbox{\halign{\hss#\hss&&\kern1.5em\hss#\hss\cr
\picture{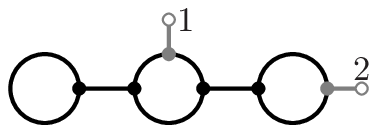}&\picture{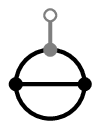}&\picture{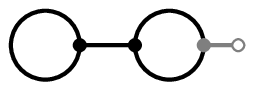}&\picture{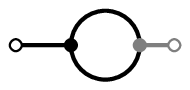}&\picture{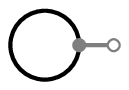}\cr
\noalign{\medskip}
(c)&(d)&(e)&(f)&(g)\cr
\crcr}}}
\figure\label{fig.e8-r}
Skeletons of reducible curves~$\B$
\endfigure
\endinsert

\subsection{Proof of Theorem~\ref{th.e8-r}}\label{proof.e8-r}
The proof repeats, almost literally, the proof of
Theorem~\ref{th.e8}. Since the two exceptional curves given by
Lemmas~\ref{val.1} and~\ref{val.2} are irreducible,
it suffices to consider a skeleton~$\Sk$ obtained
from a graph~$\Sk'$
as in
Lemma~\ref{Sk'}. Furthermore, unless $\Sk'$ is a single circle (see
Figure~\ref{fig.e8-r}\fref(g), the exceptional case in
Lemma~\ref{Sk'}), Lemma~\ref{prop.reducible} and
Corollary~\ref{2.markings} imply that $\Sk'$ is the skeleton of a
\emph{reducible} curve $B'\subset\Sigma_2$ and either $B'$ splits
into three linear components (Figure~\ref{fig.e8-r}\fref(b)
and~\fref(d)), and then
the insertion can be attached arbitrarily, or $\Sk'$ has exactly
one splitting marking
(Figure~\ref{fig.e8-r}\fref(a), \fref(c), \fref(e), and~\fref(f))
and, with respect to this marking, the insertion is attached to an
edge of type $[2,3]$ and is oriented as shown in
Figure~\ref{fig.ins}\fref(a).

\midinsert
\table\label{tab.e8-r}
Maximal sets of singularities with a type~$\bE_8$ point
represented by reducible sextics
\endtable

\def\no{}
\let\-\0
\def\*{\llap{$^*$}}
\let\*\relax
\def\tref#1{\,\text{\ref{#1}}\,}
\def\same#1{\llap{$^{#1\,}$}}
\centerline{\vbox{\offinterlineskip\halign{%
\tabstrut\ \hss#\rlap{$'$\hss}\hss\ \vrule&
\quad$#$\hss\quad\vrule&
\quad\null#\hss\quad\vrule&&\ \ \hss$#$\hss\ \ \vrule\cr
\noalign{\hrule}
\exstrut&&&&&\cr
\#&\text{Set of singularities}&
 \hss Figure&\text{Count}&
 %\text{Real?}&
 \pi_1&(l,m,n)\cr
\exstrut&&&&&\cr
\noalign{\hrule}
\exstrut&&&&&\cr
1&\bE_8\splus\bA_5\splus2\bA_3&
 \rfragment(a)1&(1,0)&\tref{s.reducible}&(4,4,6)\cr
2&\bE_8\splus\bA_7\splus\bA_3\splus\bA_1&
 \rbifrag(a)2&(0,1)&\tref{s.reducible}&(8,4,2)\cr
3&\bE_8\splus\bA_7\splus\bA_2\splus2\bA_1&
 \rfragment(a)3&(1,0)&\tref{s.leaf-out}&\cr
4&\bE_8\splus\bA_5\splus\bA_4\splus2\bA_1&
 \rfragment(b)1&(1,0)&\tref{s.reducible}&(6,5,2)\cr
5&\bE_8\splus\bA_5\splus\bA_3\splus\bA_2\splus\bA_1&
 \rfragment(b)2&(1,0)&\tref{s.reducible}&(6,3,4)\cr
6&\bE_8\splus\bA_4\splus2\bA_3\splus\bA_1&
 \rfragment(b)3&(1,0)&\tref{s.reducible}&(4,5,4)\cr
7&\same1 \bE_8\splus\bA_9\splus\bA_2&
 \rfragment(c)1&(1,0)&\tref{s.reducible}&(10,3,\-)\cr
8&\bE_8\splus\bA_9\splus2\bA_1&
 \rfragment(c)2&(1,0)&\tref{s.leaf-out}&\cr
9&\bE_8\splus\bD_8\splus\bA_2\splus\bA_1&
 \ref{fig.e8-r}\fref(d)&(1,0)&\tref{s.reducible}&(\-,3,2)\cr
10&\bE_8\splus\bD_7\splus\bA_3\splus\bA_1&
 \ref{fig.e8-r}\fref(d)&(1,0)&\tref{s.reducible}&(4,\-,2)\cr
11&\bE_8\splus\bD_6\splus\bA_3\splus\bA_2&
 \ref{fig.e8-r}\fref(d)&(1,0)&\tref{s.120}&\*(4,3,\-)\cr
12&\bE_8\splus\bD_{10}\splus\bA_1&
 \ref{fig.e8-r}\fref(e)&(1,0)&\tref{s.leaf-out}&\cr
13&\bE_8\splus\bD_6\splus\bA_5&
 \ref{fig.e8-r}\fref(e)&(1,0)&\tref{s.e8+d6+a5}&\*(6,\-,6)\cr
14&\bE_8\splus\bD_5\splus\bA_5\splus\bA_1&
 \ref{fig.e8-r}\fref(e)&(1,0)&\tref{s.leaf-out}&\cr
15&\bE_8\splus\bE_7\splus\bA_3\splus\bA_1&
 \ref{fig.e8-r}\fref(f)&(1,0)&\tref{s.leaf-out}&\cr
16&\bE_8\splus\bD_6\splus\bD_5&
 \ref{fig.e8-r}\fref(g)&(1,0)&\tref{s.e8+d6+d5}&\*\cr
17&\bE_8\splus\bE_7\splus\bD_4&
 isotrivial\!\!\!&(1,0)&\tref{s.isotrivial}&\*\cr
\exstrut&&&&&\cr
\noalign{\hrule}
\crcr}}}
\endinsert

The classification of reducible
stable maximal trigonal curves in~$\Sigma_2$
is found in~\cite{symmetric}; their skeletons
are shown (in black) in Figure~\ref{fig.e8-r}. The possible
positions of the insertion (up to symmetries of~$\Sk'$)
are shown in grey, and the corresponding sets of singularities
(see Lemma~\ref{type.E} and Remark~\ref{rem.AD})
are
listed in Table~\ref{tab.e8-r}. (The notation in the figure and
the columns in the table are the same as in the previous section.)
It turns out that each set of singularities is obtained from a
unique skeleton. (The set of singularities
$\bE_8\splus\bA_9\splus\bA_2$ prefixed with $^1$ in the table is
also realized by an irreducible sextic, see Table~\ref{tab.e8}.)
Adding the (unique) isotrivial trigonal model realizing the set of
singularities $\bE_8\splus\bE_7\splus\bD_4$, see
Proposition~\ref{1-1.e8}, and summarizing, one obtains $17$ sets
of singularities realized by $18$ curves ($16$ real ones and one
pair of complex conjugate ones).
\qed

\section{The fundamental group\label{S.groups}}

In this section, we compute the fundamental group
$\pi_1:=\pi_1(\Cp2\sminus B)$
of
a maximal sextic~$B$
%with a type~$\bE_8$ singular point.
satisfying~\cstar.
Most of the time, we assume that
the trigonal model~$\B$ of~$B$ is not isotrivial; then, we denote
by~$\GB$ the skeleton of~$\B$ and keep the notation~$u$ and~$v$
for, respectively,
the \white-- and \black-vertices spanning the insertion,
see~\ref{s.model}.

\subsection{The presentation}\label{s.presentation}
Assume that $v$ is
trivalent, choose the marking at~$v$ so that the edge $[v,u]$
is~$e_2$ at~$v$, and let $\{\Ga_1,\Ga_2,\Ga_3\}$ be a
canonical
basis in the fiber~$F_v$ over~$v$
defined by this marking (Figure~\ref{fig.basis};
see~\cite{degt.kplets} or~\cite{dessin.e7} for details).
According to~\cite{dessin.e7}, the basis
elements are subject to the following relations (the so called
\emph{relations at infinity}):
$$
\gather
\Gr^3=\Ga_1\Ga_2^2,\eqtag\label{infty.1}\\\allowdisplaybreak
\Ga_3=\Ga_2\Ga_1\Ga_2\1,\eqtag\label{infty.2}\\\allowdisplaybreak
[\Ga_1,\Ga_2^3]=1,
\eqtag\label{infty.3}
\endgather
$$
where $\Gr=\Ga_1\Ga_2\Ga_3$.
In particular, it follows that $\Ga_2^3$ is a central element.

\midinsert
\centerline{\picture{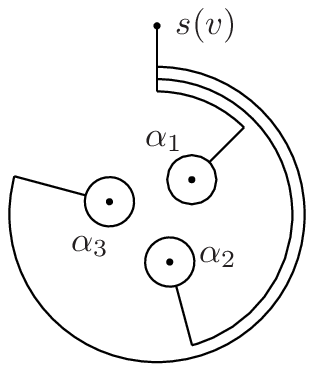}}
\figure
The canonical basis
\endfigure\label{fig.basis}
\endinsert

Note that one can eliminate~$\Ga_3$ and simplify~\eqref{infty.1}
to the following relation:
$$
(\Ga_1\Ga_2\1)^5\Ga_2^6=1.
\eqtag\label{infty.1'}
$$
We will use~\eqref{infty.1'} in the final presentations of the
groups.

Let $F_1,\ldots,F_r$ be the singular fibers of~$\B$ other
than~$F$, and let $m_j\subset\BG3$ be the braid
monodromy about~$F_j$, $j=1,\ldots,r$
(see~\cite{dessin.e7} for the choice of the reference
section and other details). Then, the Zariski--van Kampen
theorem~\cite{vanKampen} states that
%the group $\pi_1:=\pi_1(\Cp2\sminus B)$
$\pi_1$
has the following presentation
$$
\pi_1=\bigl<\Ga_1,\Ga_2,\Ga_3\bigm|
 \text{$m_j=\id$, $j=1,\ldots,r$, and
 \eqref{infty.1}--\eqref{infty.3}}\bigr>,
\eqtag\label{eq.vanKampen}
$$
where each \emph{braid relation} $m_j=\id$, $j=1,\ldots,r$,
should be understood as
the triple of relations $m_j(\Ga_i)=\Ga_i$, $i=1,2,3$.
Furthermore, in the presence of the relations at infinity, (any)
one of the braid relations $m_j=\id$, $j=1,\ldots,r$, can be
omitted.

The braid monodromy~$m_j$ can easily be computed using the
skeleton (or dessin) of~$\B$, see~\cite{degt.kplets}
or~\cite{dessin.e7} for an exposition
more tailored for the problem in
question. We omit all details and merely indicate the results.

\subsection{The monodromy to be considered}\label{s.lmn}
Given two
elements~$\Ga$, $\Gb$ of a group and a nonnegative integer~$m$,
introduce the notation
$$
\{\Ga,\Gb\}_m=\cases
(\Ga\Gb)^k(\Gb\Ga)^{-k},&\text{if $m=2k$ is even},\\
\bigl((\Ga\Gb)^k\Ga\bigr)\bigl((\Gb\Ga)^k\Gb\bigr)\1,&
 \text{if $m=2k+1$ is odd}.
\endcases
$$
The relation $\{\Ga,\Gb\}_m=1$ is equivalent to $\Gs^m=\id$, where
$\Gs$ is the generator of the braid group~$\BG2$ acting on the
free group $\<\Ga,\Gb\>$. Hence,
$$
\{\Ga,\Gb\}_m=\{\Ga,\Gb\}_n=1\quad
\text{is equivalent to}\quad
\{\Ga,\Gb\}_{\gcd(m,n)}=1.
\eqtag\label{eq.equiv}
$$
For the small values of~$m$, the relation $\{\Ga,\Gb\}_m=1$ takes
the following form:
\Dashes
\dash
$m=0$: tautology;
\dash
$m=1$: the identification $\Ga=\Gb$;
\dash
$m=2$: the commutativity relation $[\Ga,\Gb]=1$;
\dash
$m=3$: the braid relation $\Ga\Gb\Ga=\Gb\Ga\Gb$.
\endDashes

\midinsert
\centerline{\picture{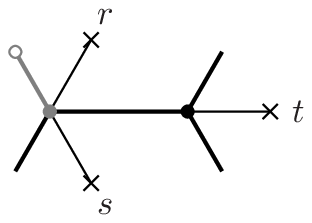}}
\figure
The braid monodromy to be considered
\endfigure\label{fig.regions}
\endinsert

Recall that the skeleton~$\GB$ of a maximal trigonal curve~$\B$ is
introduced in~\cite{degt.kplets} as a simplified version of its
\emph{dessin}~$\Gamma$, which is defined as the pull-back under
$\jB$ of the real part $\Rp1\subset\Cp1$, with the three special
values~$0$, $1$, and~$\infty$ taken into account. Thus, in
addition to \black-- and \white-vertices, $\Gamma$ also has
\cross-vertices (the pull-backs of~$\infty$); they are connected
to the \black-- and \white-vertices by, respectively, solid and
dotted edges (whereas the edges of~$\Sk$ are regarded bold). If
$\B$ is maximal, then $\Gamma$ is uniquely recovered from~$\GB$:
there is a single \cross-vertex~$a_R$ inside each
region~$R$ of~$\Sk$; it is connected to all vertices
in~$\partial R$ in the star like fashion. This vertex~$a_R$
corresponds to the unique singular fiber of~$\B$ inside~$R$, \cf\.
Remark~\ref{rem.AD}. The braid monodromy about~$a_R$ is described
in~\cite{dessin.e7}.

In most cases,
in order to compute the fundamental group, it
suffices to consider the braid monodromy about the three
\cross-vertices $r$, $s$, $t$ shown in Figure~\ref{fig.regions}.
(Here, $r$ is the \cross-vertex immediately adjacent to~$v$, $s$ is
the vertex `opposite' to~$v$, and $t$ has an extra bold edge in
the path connecting it to~$v$. Note that we do \emph{not} assume
that all three vertices are distinct.)
Assume that $r$, $s$, and~$t$ are the centers of, respectively,
an $l$-, $m$-, and~$n$-gonal region of the skeleton. Then the
relations resulting from the braid monodromy about these vertices
are
$$
\aligned
r:&\quad\{\Ga_1,\Ga_2\}_l=1,\\
s:&\quad\{\Ga_1,\Ga_s\}_m=1,\quad
 \text{where $\Ga_s=\Ga_2\Ga_3\Ga_2\1$},\\
t:&\quad\{\Ga_2,\Ga_t\}_n=1,\quad
 \text{where $\Ga_t=(\Ga_1\Ga_2)\Ga_3(\Ga_1\Ga_2)\1$},
\endaligned\eqtag\label{eq.lmn}
$$
{\em provided that the singular
fiber of~$\B$ over the corresponding vertex
is of type~$\tA{}$}. (If the fiber is of type~$\tD{}$, we usually
ignore the corresponding relation.)
The additional relations~\eqref{eq.lmn}, as well as the original
relations at infinity \eqref{infty.1}--\eqref{infty.3}, are
easily programmable in \GAP~\cite{GAP}, see Figure~\ref{GAP.in}.

\midinsert
\beginGAP%
gap> a := FreeGroup(3);;
gap> ## Common relations (at infinity)
gap> r1 := (a.1*a.2*a.3)^3/(a.1*a.2^2);;
gap> r2 := a.2*a.1*a.2^-1/a.3;;
gap> r3 := a.1*a.2^3/(a.2^3*a.1);;
gap> ## Common elements to be used in other relations
gap> as := a.2*a.3*a.2^-1;;
gap> at := (a.1*a.2)*a.3*(a.1*a.2)^-1;;
gap> ## The variable relation (n=0 means no relation)
gap> r := function(a,b,n)
> if (n mod 2) = 0 then
>  return (a*b)^(n/2)/(b*a)^(n/2);
> else
>  return (a*b)^((n-1)/2)*a/((b*a)^((n-1)/2)*b);
> fi;
> end;;
gap> ## Most fundamental groups are handled by this
gap> size := function(l, m, n)
> return Size(a/[r1, r2, r3, r(a.1,a.2,l), r(a.1,as,m), r(a.2,at,n)]);
> end;;
gap> size2 := function(l, m, n)
> return Size(a/[r1, r2, r3, r(a.1,a.2,l), r(a.1,as,m), r(a.2,at,n), a.2^3]);
> end;;
\endGAP
\figure
The \GAP\ input
\endfigure\label{GAP.in}
\endinsert

Note that any of relations~\eqref{eq.lmn} can easily be ignored:
one should just let
the corresponding parameter~$l$, $m$, or~$n$ to be equal
to zero. To emphasize the fact
that a relation is omitted, we will use `$\0$' instead
of~`$0$' in the references to~\eqref{eq.lmn}.

\paragraph\label{s.others}
For most irreducible maximal curves, the group~$\pi_1$ is computed
using \GAP~\cite{GAP}: the function {\tt size(l,m,n)} in
Figure~\ref{GAP.in} returns~{\tt6} and, since
$\pi_1/[\pi_1,\pi_1]=\CG6$,
this implies that $\pi_1$ is abelian. The values $(l,m,n)$
used are listed in Table~\ref{tab.e8}; they are easily read from
the skeletons of the curves, see Figure~\ref{fig.e8}.
The $\tD{}$-type singular fibers, if
present, are ignored.

\paragraph\label{s.reducible}
For reducible maximal curves, we use the following obvious
observation: {\proclaimfont let~$G$ be a group, and let
$H\subset G$ be a central subgroup such that the projection
$H\to G/[G,G]$ is a monomorphism. Then the commutants of~$G$ and
of~$G/H$ are isomorphic.} We apply this statement to the subgroup
$H\subset\pi_1$ generated by the central element~$\Ga_2^3$ and
compute the quotient $\pi_1/\Ga_2^3$ using \GAP~\cite{GAP}. Note
that, if the curve is known to be reducible, than $\Ga_2^3$
projects to an element of infinite order in
$\pi_1/[\pi_1,\pi_1]$ (hence, the statement above does apply),
%the abelianization
%of~$\pi_1$
and the abelianization of $\pi_1/\Ga_2^3$ is $\CG{15}$. Thus, if
the function {\tt size2(l,m,n)} in Figure~\ref{GAP.in} returns
{\tt15}, we conclude that $\pi_1$ is abelian. The values $(l,m,n)$
used are listed in Table~\ref{tab.e8-r}; they are read from the
skeletons of the curves, see Figure~\ref{fig.e8-r},
with the $\tD{}$-type singular fibers
ignored.

In the sequel, without further references, we assume
that all statements like `$\ord\pi_1=6$' or
`$\ord(\pi_1/\Ga_2^3)=15$' are proved using~\GAP~\cite{GAP}.

\subsection{Insertion close to a loop}
Here, we consider a few special positions of the insertion with
respect to a loop (\ie, a monogonal region)
of~$\Sk'$.

\midinsert
\centerline{\vbox{\halign{\hss#\hss&&\kern6em\hss#\hss\cr
\cpic{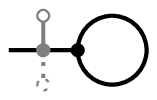}&\cpic{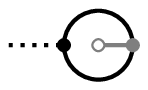}&\cpic{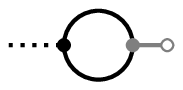}\cr
\noalign{\medskip}
(a)&(b)&(c)\cr
\crcr}}}
\figure
Fragments resulting in abelian~$\pi_1$
\endfigure\label{fig.fragments}
\endinsert

\paragraph\label{s.stem}
Assume that the skeleton of~$\B$ has a fragment shown in
Figure~\ref{fig.fragments}\fref(a), \ie, the insertion is right next to
a loop.
Then $\pi_1$ has
relation~\eqref{eq.lmn} with $(l,m,n)=(\0,\0,1)$,
and
%using \GAP~\cite{GAP}
one concludes that $\pi_1=\CG6$
%, see Figure~\ref{GAP.out}.
(as {\tt size(0,0,1)} returns {\tt 6}).
If the
insertion is as shown in Figure~\ref{fig.fragments}\fref(a) by the
dotted lines, then still $\pi_1=\CG6$: one can use either a
similar calculation or just the symmetry arguments.
%In particular,
As a consequence,
the curve~$B$ is necessarily irreducible in this case.

\paragraph\label{s.leaf-in}
Assume that the skeleton of~$\B$ has a fragment shown in
Figure~\ref{fig.fragments}\fref(b), \ie, the insertion is inside a
loop of~$\Sk'$.
(The rightmost \black-vertex in the figure
can be either bi- or trivalent; it is
not used.)
%(The dotted bold line indicates
%that the valency of the \black-vertex is irrelevant: it can be
%either~$2$ or~$3$.)
Then $\pi_1$ has
relation~\eqref{eq.lmn} with $(l,m,n)=(3,\0,\0)$,
and
%using \GAP~\cite{GAP}
one concludes that $\pi_1=\CG6$
(as {\tt size(3,0,0)} returns {\tt 6}). This fact implies, in
particular, that the curve~$B$ is irreducible.

%\subsection{Insertion outside a loop}\label{s.leaf-out}
\paragraph\label{s.leaf-out}
Assume that the skeleton of~$\B$ has a fragment shown in
Figure~\ref{fig.fragments}\fref(c), \ie, the insertion is
right outside a loop of~$\Sk'$.
(The rightmost \black-vertex in the figure
can be either bi- or trivalent; it is
not used.)
%(The dotted bold line indicates
%that the valency of the \black-vertex is irrelevant: it can be
%either~$2$ or~$3$.)
Then the group $\pi_1$ has
relation~\eqref{eq.lmn} with $(l,m,n)=(\0,2,\0)$, \ie,
$[\Ga_1,\Ga_2\Ga_3\Ga_2\1]=1$. Using~\eqref{infty.2}
and~\eqref{infty.3}, this relation can be rewritten in the form
$[\Ga_1,\Ga_2\1\Ga_1\Ga_2]=1$. Hence, also
$[\Ga_1,\Ga_2\Ga_1\Ga_2\1]=1$. Spelling~\eqref{infty.1} out and
eliminating~$\Ga_3$ with the help of~\eqref{infty.2},
after the cancellation
one has
$$
\Ga_1\cdot\Ga_2\1\Ga_1\Ga_2\cdot\Ga_2\Ga_1\Ga_2\1\cdot
\Ga_1\cdot\underline{\Ga_2}\cdot\Ga_2\Ga_1\Ga_2\1=1.
$$
Thus, $\Ga_2$ (the underlined instance) is a product of elements
commuting with~$\Ga_1$, and the group is abelian.

Observe that these arguments apply to both reducible and
irreducible curves.

\subsection{A few special cases}
Below, we treat the few special cases that are not covered
by~\ref{s.lmn} directly.

\subsubsection{The set of singularities $\bE_8\splus\bD_6\splus\bA_5$
\rm(No\.~$13'$ in Table~\ref{tab.e8-r})}\label{s.e8+d6+a5}
The
type~$\tD{_6}$ fiber
over~$s$ prevents one from using relation~\eqref{eq.lmn} with
$m=2$. However, the relations at
infinity~\eqref{infty.1}--\eqref{infty.3},
relations~\eqref{eq.lmn} with $(l,m,n)=(6,\0,6)$, and the
relation
$$
\Ga_3=(\Ga_2\Ga_1\Ga_2)\1\Ga_1(\Ga_2\Ga_1\Ga_2)
$$
resulting from the braid monodromy about the leftmost loop in
Figure~\ref{fig.e8-r}\fref(e) suffice to show
that $\ord(\pi_1/\Ga_2^3)=15$. Hence, $\pi_1$ is abelian,
see~\ref{s.reducible}.

%\paragraph\label{s.e8+d6+d5}
\subsubsection{The set of singularities $\bE_8\splus\bD_6\splus\bD_5$
\rm(No\.~$16'$ in Table~\ref{tab.e8-r})}\label{s.e8+d6+d5}
The $\tD{}$-type fibers
prevent one from using~\eqref{eq.lmn}.
However,
the type~$\tD{_6}$ fiber over~$r$ gives a relation
$[\Ga_3,\Ga_1\Ga_2]=1$. Then, $[\Ga_3,\Ga_2]=1$
(from~\eqref{infty.1}) and $\Ga_3=\Ga_1$ (from~\eqref{infty.2});
hence, the group is abelian.
%adding to~\eqref{infty.1}--\eqref{infty.3} the relation
%$[\Ga_3,\Ga_1\Ga_2]=1$ arising from the monodromy about the
%type~$\tD{_6}$ fiber over~$r$, one can see
%that $\ord(\pi_1/\Ga_2^3)=15$. Hence, $\pi_1$ is abelian,
%see~\ref{s.reducible}.

%\paragraph\label{s.2e8+a3}
\subsubsection{The set of singularities $2\bE_8\splus\bA_3$
\rm(No\.~29 in Table~\ref{tab.e8})}\label{s.2e8+a3}
The vertex~$v$ is bivalent, and we choose the reference fiber over
an inner point of the edge $[u,v]$. Relations at infinity
\eqref{infty.1}--\eqref{infty.3} still hold, and the
type~$\tE{_8}$ fiber over~$v$ gives, among others, the relation
$\Ga_2=\Gr^2\Ga_3\Gr^{-2}$. Now,
%\GAP~\cite{GAP} shows
one can see
that $\pi_1=\CG6$.

%\paragraph\label{s.e8+e6+d5}
\subsubsection{The set of singularities $\bE_8\splus\bE_6\splus\bD_5$
\rm(No\.~30 in Table~\ref{tab.e8})}\label{s.e8+e6+d5}
The vertex~$v$ is monovalent.
As in the previous case, choose the reference fiber over
an inner point of $[u,v]$.
Then, it suffices to add to \eqref{infty.1}--\eqref{infty.3} the
relation $\Ga_3=\Gr\Ga_2\Gr\1$ arising from the monodromy about
the type~$\tE{_6}$ fiber over~$v$
%and use \GAP~\cite{GAP}
to conclude that $\pi_1=\CG6$.

\subsection{The set of singularities
$\bE_8\splus\bA_4\splus\bA_3\splus2\bA_2$}\label{s.720}
%(No\.~1 in Table~\ref{tab.e8})}\label{s.720}
This is No\.~$1$ in Table~\ref{tab.e8}.
As explained in~\ref{s.presentation}, one of the four braid
relations corresponding to the four singular fibers other than~$F$
can be ignored. Ignoring the upper left triangle
in Figure~\ref{fig.e8}\fref(a), one
%can see that the group~$\pi_1$
%is given by
arrives at the presentation
$$
\bigl<\Ga_1,\Ga_2,\Ga_3\bigm|
\text{\eqref{infty.1}--\eqref{infty.3} and
 \eqref{eq.lmn} with $(l,m,n)=(5,4,3)$}\bigr\>.
\eqtag\label{eq.G}
$$
Denote this group by~$G_6$. Using \GAP~\cite{GAP},
see Figure~\ref{GAP.a4+a3+2a2}, one can see that:
\roster
\item\local{G.1}
one has $\ord G_6=720$, and $[G_6,G_6]$ is a perfect group
of order $120$;
\item\local{G.2}
the only perfect group of order $120=\ord[G_6,G_6]$ is
$\SL(2,\Bbb F_5)$;
\item\local{G.3}
the last two relations in~\eqref{eq.G} follow from the others (as
dropping these relations does not change the
order of the group);
\item\local{G.4}
$G_6$ is a semidirect product of $G_6/[G_6,G_6]=\CG6$ and
$[G_6,G_6]=\SL(2,\Bbb F_5)$;
\item\local{G.5}
the order of each generator $\Ga_i$, $i=1,2,3$, in~$G$ equals
$12$;
\item\local{G.6}
the group is generated by $\Ga_2\Ga_3\Ga_2\1$, $\Ga_1$, and
$\Ga_3$;
\item\local{G.7}
the centralizer~$C$ of~$[G_6,G_6]$ is isomorphic to~$\CG{12}$ and
$C\cap[G_6,G_6]=\CG2$; in particular,
the canonical projection $C\to G_6/[G_6,G_6]$ is onto.
\endroster
The presentation of~$G_6$ stated in Theorem~\ref{th.main} is
obtained from~\eqref{eq.G} by dropping the last two relations,
see statement~\loccit{G.3} above, eliminating
the last generator~$\Ga_3$ \via~\eqref{infty.2},
and replacing~\eqref{infty.1} with~\eqref{infty.1'}. The central
product description of the group given in
Remark~\ref{rem.groups} follows from
statement~\loccit{G.7} above: one has
$$
G_6=\bigl(C\times[G_6,G_6]\bigr)/\bigl(C\cap[G_6,G_6]\bigr).
$$

\midinsert
\beginGAP%
gap> ## l=5, m=4, n=3 (no.1)
gap> g := a / [r1, r2, r3, r(a.1,a.2,5), r(a.1,as,4), r(a.2,at,3)];;
gap> Size(g);
720
gap> Size(a / [r1, r2, r3, r(a.1,a.2,5)]);             ## drop two relations
720
gap> List(DerivedSeriesOfGroup(g), Size);              ## size of [G,G]
[ 720, 120 ]
gap> List(DerivedSeriesOfGroup(g), AbelianInvariants); ## [G,G] is perfect
[ [ 2, 3 ], [  ] ]
gap> NumberPerfectGroups(120); PerfectGroup(120);
1
A5 2^1
gap> StructureDescription(DerivedSubgroup(g));         ## [G,G]
"SL(2,5)"
gap> StructureDescription(g);                          ## semidirect product
"C3 x (SL(2,5) : C2)"
gap> Size(Subgroup(g, [g.1]));                         ## order of g.1
12
gap> Index(g, Subgroup(g, [g.2*g.3*g.2^-1, g.1, g.3]));## generating subset
1
gap> h := DerivedSubgroup(g);; c := Centralizer(g, h);;## centralizer of G'
gap> Size(Intersection(c, h));
2
gap> StructureDescription(c);
"C12"
\endGAP
\figure
The \GAP\ output, the set of singularities
$\bE_8\splus\bA_4\splus\bA_3\splus2\bA_2$
\endfigure\label{GAP.a4+a3+2a2}
\endinsert

\subsection{The set of singularities
$\bE_8\splus\bD_6\splus\bA_3\splus\bA_2$}\label{s.120}
This is No\.~$11'$ in Table~\ref{tab.e8-r}.
As explained in~\ref{s.presentation}, one of the three braid
relations corresponding to the three singular fibers other than~$F$
can be ignored. Ignoring the type~$\tD{_6}$ fiber over~$t$,
one arrives at the presentation
$$
\bigl<\Ga_1,\Ga_2,\Ga_3\bigm|
\text{\eqref{infty.1}--\eqref{infty.3} and
 \eqref{eq.lmn} with $(l,m,n)=(4,3,\0)$}\bigr\>.
\eqtag\label{eq.Gi}
$$
Denote this group by~$G_\infty$, and analyze
the quotient $G=G_\infty/\Ga_2^3$
using \GAP~\cite{GAP}, see Figure~\ref{GAP.3a1}. One has:
\roster
\item\local{Gi.1}
$\ord G=1800$, and $[G,G]$ is a perfect group
of order $120$;
\item\local{Gi.3}
any of the last two relations in~\eqref{eq.Gi}
follows from the other relations (as
dropping a relation does not change the
order of the group);
\item\local{Gi.6}
the group~$G$
is generated by $\Ga_2\Ga_3\Ga_2\1$, $\Ga_1$, and $\Ga_3$; hence,
$G_\infty$ is generated by $\Ga_2\Ga_3\Ga_2\1$, $\Ga_1$, $\Ga_3$,
and the central element~$\Ga_2^3$;
\item\local{Gi.7}
the centralizer~$C$ of~$[G,G]$ is isomorphic to~$\CG{30}$ and
$C\cap[G,G]=\CG2$; in particular,
the canonical projection $C\to G/[G,G]$ is onto.
\endroster
Hence, one has $[G_\infty,G_\infty]=[G,G]=\SL(2,\Bbb F_5)$,
\cf\.~\ref{s.reducible} and~\iref{s.720}{G.2}, and, since the
abelianization
$G_\infty/[G_\infty,G_\infty]=\Z$ is free, $G_\infty$ splits into
a semi-direct product.
From statement~\loccit{Gi.7} above it
follows that the generator of the abelianization lifts to an
element commuting with $[G_\infty,G_\infty]$; hence, the product
is in fact direct, see Remark~\ref{rem.groups}.
The presentation of~$G_\infty$ stated in Theorem~\ref{th.main} is
obtained from~\eqref{eq.Gi} by dropping the last relation
$\{\Ga_1,\Ga_s\}_3=1$,
see statement~\loccit{Gi.3} above, eliminating
the last generator~$\Ga_3$ \via~\eqref{infty.2},
and replacing \eqref{infty.1} with~\eqref{infty.1'}.

\midinsert
\beginGAP%
gap> g := a/[r1, r2, r3, r(a.1,a.2,4), r(a.1,as,3), a.2^3];;
gap> Size(g);
1800
gap> List(DerivedSeriesOfGroup(g), Size);              ## size of [G,G]
[ 1800, 120 ]
gap> List(DerivedSeriesOfGroup(g), AbelianInvariants); ## [G,G] is perfect
[ [ 3, 5 ], [  ] ]
gap> Size(a/[r1, r2, r3, r(a.1,a.2,4), a.2^3]);        ## drop a relation
1800
gap> Size(a/[r1, r2, r3, r(a.1,as,3), a.2^3]);         ## drop a relation
1800
gap> Index(g, Subgroup(g, [g.2*g.3*g.2^-1, g.1, g.3]));## generating subset
1
gap> h := DerivedSubgroup(g);; c := Centralizer(g, h);;## centralizer of G'
gap> Size(Intersection(c, h));
2
gap> StructureDescription(c);
"C30"
\endGAP
\figure
The \GAP\ output, the set of singularities
$\bE_8\splus\bD_6\splus\bA_3\splus\bA_2$
\endfigure\label{GAP.3a1}
\endinsert

\subsection{Isotrivial curves}\label{s.isotrivial}
The trigonal model~$\B$ of a plane sextic~$B$ with
the set of singularities $\bE_8\splus\bE_7\splus\bD_4$
(No\.~$17'$ in Table~\ref{tab.e8-r}) is isotrivial, and it is the
only isotrivial trigonal model of a maximal sextic, see
Proposition~\ref{1-1.e8}.

One has $j_{\B}\equiv1$ and, in
appropriate affine coordinates $(x,y)$ in~$\Sigma_3$, the
Weierstra{\ss} equation of~$\B$ has the form
$$
y^3+p(x)y=0,
\eqtag\label{eq.j=1}
$$
where $\deg p=5$. (We assume that the distinguished fiber~$F$,
corresponding to a simple root of~$p$, is at $x=\infty$.)

In general, for any curve~$\B$ given by~\eqref{eq.j=1}, the braid
monodromy is abelian: the singular fibers of~$\B$ are in a
one-to-one correspondence with the roots of~$p$, and the monodromy
$m_j$ about the fiber~$F_j$ corresponding to a $s_j$-fold root
of~$p$ is $(\Gs_1\Gs_2\Gs_1)^{s_j}$ (for an appropriate basis in
the reference fiber).
Hence, the braid
relations are equivalent to a single relation
$(\Gs_1\Gs_2\Gs_1)^s=\id$, where $s$ is the greatest
common divisor of the multiplicities of all roots of~$p$.

If $\B$ as above is the trigonal model of a plane sextic
(not necessarily maximal) given by
Proposition~\ref{1-1.e8}, then $\deg p=5$ (a simple root at
infinity) and the multiplicity of each root of~$p$ is at most~$3$
(simple singularities only); hence, the greatest common
divisor~$s$ above equals~$1$, and the resulting relation
$\Gs_1\Gs_2\Gs_1=\id$ yields $\Ga_1=\Ga_3$, $[\Ga_1,\Ga_2]=1$.
Thus, the fundamental group is abelian.

\section{Perturbations}\label{S.perturbations}

We start with a description of the fundamental groups of the
perturbations of a singular point of type~$\bE_8$ or~$\bD_m$. We
apply these results to
%the study of
the perturbations of
maximal sextics satisfying~\cstar, proving
Theorems~\ref{th.main} and~\ref{th.pert}.

\subsection{Perturbations of a type~$\bE_8$ singular point}\label{s.pert.e8}
Consider a type~$\bE_8$ singular point~$P$
of a plane curve~$B$
and let~$\MB$ be a Milnor ball about~$P$. Let $B_t$ be a
perturbation of $B=B_0$ transversal to the boundary~$\partial\MB$.
We are interested in the perturbation epimorphism
$\pi_1(\MB\sminus B)\onto\pi_1(\MB\sminus B_1)$.

A local normal form of~$B$ at~$P$ is $\{y^3+x^5=0\}$. Consider the
line $L_\Ge=\{x=\Ge^3\}$, where $\Ge>0$ is a real number,
$\Ge^3\ll\text{radius of~$\MB$}$. The intersection $L_\Ge\cap B$
consists of the vertices $y_k=-\Ge^5\exp(2\pi ki/3)$, $k=0,1,2$,
of an equilateral triangle, and we denote by
$\{\bc_1,\bc_2,\bc_3\}$ a corresponding canonical basis for the
group
$\pi_1(L_\Ge\sminus B)$ (\cf\. Figure~\ref{fig.basis}); clearly,
$\{\bc_1,\bc_2,\bc_3\}$ is
also a basis for $\pi_1(\MB\sminus B)$, and one has
$$
\pi_1(\MB\sminus B)=\bigl<\bc_1,\bc_2,\bc_3\bigm|
 \bc_1\Gr^2=\Gr^2\bc_2,\
 \bc_2\Gr^2=\Gr^2\bc_3,\
 \bc_3\Gr=\Gr\bc_1\bigr>,
$$
where $\Gr=\bc_1\bc_2\bc_3$.

\proposition\label{pert.e8}
Up to deformation, there are
three proper perturbations of a type~$\bE_8$
singularity with nonabelian fundamental group.
They are as follows\rom:
\Dashes
\dash
$\bA_4\splus\bA_3$\rom:
$\{\bc_1,\bc_2\}_4=\{\bc_1,\bc_3\}_5=1$ and
$\bc_1=\bc_3\bc_1\bc_2\bc_1\1\bc_3\1$\rom;
\dash
$\bA_4\splus\bA_2\splus\bA_1$\rom:
$\{\bc_1,\bc_2\}_5=\{\bc_1,\bc_3\}_3=[\bc_1\bc_2\bc_1\1,\bc_3]=1$\rom;
\dash
$\bD_5\splus\bA_2$\rom:
$\{\bc_1,\bc_2\}_3=[\bc_2,\bc_3\bc_1]=1$ and
$\bc_3=(\bc_1\bc_2\bc_3)\bc_1(\bc_1\bc_2\bc_3)\1$,
\endDashes
where listed are the sets of singularities of the perturbed
curves~$B_t$, $t\ne0$, and the relations in the group
$\pi_1(\MB\sminus B_t)$ in the basis $\{\bc_1,\bc_2,\bc_3\}$
described above.
\endproposition

\proof
According to E.~Looijenga~\cite{Looijenga}, the deformation
classes of perturbations of a simple singularity are in a
one-to-one correspondence with the induced subgraphs of its Dynkin
graph. In particular, there are eight maximal (\ie, those with the
total Milnor number $\mu=7$) perturbations of~$\bE_8$. We assert
that each of these eight perturbations
can be realized by a maximal trigonal curve
$\B\subset\Sigma_2$ with a type~$\tA{_0^{**}}$ singular fiber~$F$
(which we place at infinity and cut off) and all other
fibers stable.

\midinsert
\centerline{\vbox{\halign{\hss#\hss&&\qquad\quad\hss#\hss\cr
\cpic{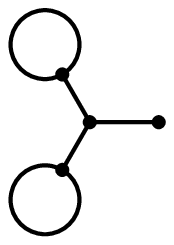}&\cpic{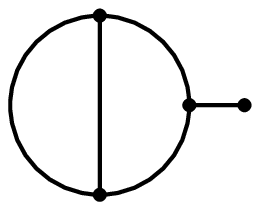}&\cpic{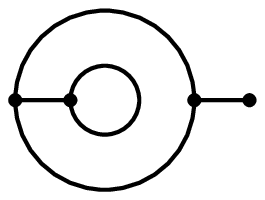}\cr
\noalign{\medskip}
$\bA_7$&$\bA_4\splus\bA_2\splus\bA_1$&$\bA_4\splus\bA_3$\cr
\crcr}}}
\bigskip
\centerline{\vbox{\halign{\hss#\hss&&\qquad\quad\hss#\hss\cr
\cpic{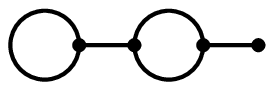}&\cpic{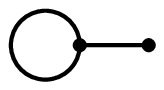}&\cpic{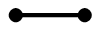}&\cpic{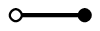}\cr
\noalign{\medskip}
$\bA_6\splus\bA_1$&$\bD_7$ and $\bD_5\splus\bA_2$&$\bE_6\splus\bA_1$&$\bE_7$\cr
\crcr}}}
\figure\label{fig.p-e8}
The perturbations of~$\bE_8$
\endfigure
\endinsert

Indeed, there are eight such maximal curves,
defined by the seven skeletons shown in Figure~\ref{fig.p-e8}, and
their sets of singularities are exactly those predicted
by~\cite{Looijenga}. (The skeletons in Figure~\ref{fig.p-e8} are
characterized by the fact that they have a monovalent
\black-vertex; formally,
any such skeleton can be obtained by shrinking a loop
in the skeleton of a stable maximal curve in~$\Sigma_2$,
see~\cite{symmetric}; this observation gives one a complete
classification.)
On the other hand, the Weierstra{\ss}
equation of any trigonal curve $\B\subset\Sigma_2$
with a type~$\tA{_0^{**}}$ singular fiber~$F$ at infinity
has the form
$$
f(x,y):=y^3+p(x)y+q(x)=0,
$$
where $\deg p=3$ and $\deg q=5$ (\ie, both~$p$ and~$q$ have a
simple root at infinity). This equation can be renormalized to
make the leading coefficient of~$q$ equal to one, and then the
family~$\B_t$ given by the polynomial
$f_t(x,y)=t^{15}f(x/t^3,y/t^5)$ defines a degeneration of
$\B=\B_1$ to the isotrivial curve $\B_0=\{y^3+x^5=0\}$ with a type
$\bE_8$ singular point at the origin. (The family is indeed a
degeneration as the curves~$\B_1$ and~$\B_t$, $t\ne0$, are related
by the automorphism $(x,y)\mapsto(t^3x,t^5y)$ of~$\Sigma_2$; note
that this automorphism preserves the
type~$\tA{_0^{**}}$ fiber~$F$ at infinity.)

Pick a Milnor ball~$\bar\MB$ about the type
$\bE_8$ singular point of~$\B$ and consider the pair
$(\bar\MB,\B)$ instead of $(\MB,B)$. Since $\B$ has no singular
fibers except~$0$ and~$\infty$, there is a diffeomorphism
$\bar\MB\sminus\B\cong\Sigma_2\sminus(\B\cup E\cup F)$. Then,
since the perturbation $\B_t$ constructed above is `constant' in a
neighborhood of infinity, there also are diffeomorphisms
$\bar\MB\sminus\B_t\cong\Sigma_2\sminus(\B_t\cup E\cup F)$. (We
assume that all curves~$\B_t$ are transversal to
$\partial\bar\MB$.)

Thus, it remains to compute the groups
$\pi_1(\Sigma_2\sminus(\B_1\cup E\cup F))$. The computation is
very similar to~\S\ref{S.groups}; it uses Zariski--van Kampen's
approach~\cite{vanKampen}.
Pick a trivalent \black-vertex~$v$ of~$\Sk$.
(In most cases, we take
for~$v$ the vertex adjacent to
the distinguished monovalent \black-vertex~$u$
corresponding to~$F$.
The two exceptional cases, when $\Sk$ has no trivalent
\black-vertices, are treated separately.)
%monovalent \black-vertex~$u$ of~$\Sk$ is adjacent to a trivalent
%\black-vertex~$v$.
%(The two exceptional cases are treated below
%separately.)
Let $F_v$ be the fiber over~$v$, and let
$\{\bb_1,\bb_2,\bb_3\}$ be a canonical basis for the group
$\pi_1(F_v\sminus(\B_1\cup E))$
defined by the marking
with respect to which $[u,v]$ is the edge~$e_1$ at~$v$, \cf\.
Figure~\ref{fig.basis}. Then the group
$\pi_1(\Sigma_2\sminus(\B_1\cup E\cup F))$
has a presentation
$$
\pi_1(\Sigma_2\sminus(\B_1\cup E\cup F))=
 \<\bb_1,\bb_2,\bb_3\,|\,\text{$m_j=\id$, $j=1,\ldots$}\>,
$$
and the braid monodromies~$m_j$ about the singular fibers
of~$\B_1$ are computed as explained in~\cite{dessin.e7}
or~\cite{degt.kplets}. Note that the presentation above has no
relation at infinity; this is due to the fact that the fiber~$F$
at infinity remains removed.

Below, we consider the eight maximal perturbations one by one. By
default, $v$ is the vertex adjacent to~$u$. To shorten the
notation,
we abbreviate the group
$\pi_1(\Sigma_2\sminus(\B_1\cup E\cup F))$ in question by~$G$.

\subsubsection{The perturbation~$\bA_7$}
Take for~$v$ the corner of the upper monogonal region in
Figure~\ref{fig.p-e8}.
Then the relations are $\bb_1=\bb_3$,
$\bb_2=\bb_3\1\bb_2\1\bb_1\bb_2\bb_3$, and $\{\bb_1,\bb_2\}_8=1$.
Eliminating~$\bb_3$, one can simplify the second relation to
$\{\bb_1,\bb_2\}_3=1$. Hence, $G$ is abelian due
to~\eqref{eq.equiv}.

\subsubsection{The perturbation $\bA_6\splus\bA_1$}
The relations are
$[\bb_2,\bb_3]=\{\bb_1,\bb_2\}_7=1$ and $\bb_3=\bb_1\bb_2\bb_1\1$.
Thus, $G$ is generated by~$\bb_1$, $\bb_2$, and the second
relation implies that $(\bb_1\bb_2)^7$ is a central element. Using
\GAP~\cite{GAP}, one can see that
$\ord(G/(\bb_1\bb_2)^7)=14$; hence, this quotient is abelian, and
so is~$G$, \cf\.~\ref{s.reducible}.

\subsubsection{The perturbation~$\bD_7$}
Among other relations, one has
$\bb_2=\bb_3$ (from the vertical tangent) and
$[\bb_3,\bb_1\bb_2]=1$ (from the type~$\tD{_7}$ fiber).
Eliminating~$\bb_3$, one concludes that $G$ is abelian.

\subsubsection{The perturbation $\bE_6\splus\bA_1$}
The skeleton~$\Sk$ has
no trivalent \black-vertices, and we choose the reference
fiber~$F_v$ over a point~$v$ in the solid edge of the dessin of
the curve~$\B_1$
connecting its
type~$\tE{_6}$ and type~$\tA{_1}$ singular fibers. Let
$\{\bb_1,\bb_2,\bb_3\}$ be an appropriate `canonical' basis
in~$F_v$, such that the generators~$\bb_2$ and~$\bb_3$ are brought
together when the fiber approaches the node. Then the relations
are $[\bb_2,\bb_3]=1$ (from the node),
$\bb_2=\rho\bb_1\rho\1$, and
$\bb_3=\rho\bb_2\rho\1$, where $\rho=\bb_1\bb_2\bb_3$. Conjugating
the first relation by~$\rho$ and taking into account the last two
relations, one obtains $[\bb_1,\bb_2]=1$. Hence, $\bb_2$ is a
central element,
%$[\bb_2,\rho]=1$,
and the last relation implies
$\bb_2=\bb_3$. Thus, $G$ is abelian.

\subsubsection{The perturbation~$\bE_7$}
This time, we choose the reference
fiber~$F_v$ over a point~$v$ in the dotted edge of the dessin of
the curve~$\B_1$
connecting its
type~$\tE{_7}$ and type~$\tA{_0^*}$ singular fibers, and take for
$\{\bb_1,\bb_2,\bb_3\}$ an appropriate `linear'
basis in~$F_v$, so that $\bb_1$ and~$\bb_2$ are brought together
when the point approaches the vertical tangent~$\tA{_0^*}$.
Among other
relations, one has $\bb_1=\bb_2$ (from the vertical tangent)
and
$[\bb_2,\bb_1\bb_2\bb_3\bb_1]=1$. Eliminating~$\bb_2$, one
concludes that $G$ is abelian.

\subsubsection{Other maximal perturbations}
For the remaining three maximal perturbations, we merely list the
relations for~$G$. They are as follows:
\Dashes
\dash
$\bA_4\splus\bA_3$:
$\{\bb_1,\bb_2\}_4=\{\bb_2,\bb_3\}_5=1$ and
$\bb_2=\bb_3\bb_1\bb_3\1$;
%$\ord(G/\bb_1^6)=11520$;
\dash
$\bA_4\splus\bA_2\splus\bA_1$:
$\{\bb_1,\bb_2\}_5=\{\bb_2,\bb_3\}_3=[\bb_1,\bb_3]=1$;
%$\ord(G/\bb_1^2)=120$;
\dash
$\bD_5\splus\bA_2$:
$\{\bb_1,\bb_2\}_3=[\bb_1,\bb_2\bb_3]=1$ and
$\bb_3=(\bb_1\bb_2\bb_3)\bb_2(\bb_1\bb_2\bb_3)\1$.
%$\ord(G/\bb_1^5)=600$.
\endDashes
All three groups are nonabelian: in the order of appearance, one
has
$\ord(G/\bb_1^3)=360$, $\ord(G/\bb_1^2)=120$, and
$\ord(G/\bb_1^5)=600$. (If $G$ were abelian,
for any integer $n>0$ one would have
$G/\bb_1^n=\CG{n}$.)
Note that, for the last group (the
perturbation $\bD_5\splus\bA_2$), one also has
$$
\ord(G/\bb_1^{12})=12,\quad\text{hence}\quad
G/\bb_1^{12}=\CG{12}.
\eqtag\label{eq.G12}
$$
To complete the proof for the maximal perturbations, it remains to
notice that, from the point of view of the trigonal implementation
of the type~$\bE_8$ singularity, the line~$L_\Ge$ introduced at
the beginning of this section (the line carrying the basis
$\{\bc_1,\bc_2,\bc_3\}$) should be regarded as a fiber `close to
infinity', \eg, the fiber over a point in the edge $[u,v]$
of~$\Sk$ close
to~$u$. It is related to the fiber~$F_v$ used in the computation
\via\ the monodromy through the \white-vertex at the middle of the
edge; hence, the two bases are related \via\
$\bb_1=\bc_1\bc_2\bc_1\1$,
$\bb_2=\bc_1$, $\bb_3=\bc_3$. Substituting, one obtains the
presentations announced in the statement.

\subsubsection{Non-maximal perturbations}
We assert that the fundamental group of any perturbation with the
total Milnor number $\mu=6$ (and hence of any non-maximal
perturbation) is abelian. For proof, one can list all such
perturbations (obtained by removing two vertices from the Dynkin
graph~$\bE_8$) and show, \eg, using trigonal curves and their
dessins, that each of them degenerates to a maximal one with
abelian group. Alternatively, if a perturbation appears to
degenerate to a maximal one with nonabelian group, one can analyze
the extra relations (using the dessins again) and show that the
perturbed group is abelian. We omit the details.
\endproof

\subsection{Perturbations of $\bD$-type singular points}
Consider a type~$\bD_m$, $m>4$, singular point~$Q$ of a plane
curve~$B$. According to Looijenga~\cite{Looijenga}, its
perturbations are classified by the induced subgraphs of the
Dynkin graph~$\bD_m$. In particular, maximal are the perturbation
$\bA_{m-1}$ (removing a short end of the diagram) or $\bD_{m-1}$
and $\bD_p\splus\bA_{m-p-1}$, $2\le p\le m-2$ (removing a vertex
from the long end or the trivalent vertex); in the latter case, in
order to emphasize the perturbation, we let $\bD_2=2\bA_1$ and
$\bD_3=\bA_3$.

All perturbations can be realized by trigonal curves. Thus, we
assume that $Q$ is a singular point of a trigonal curve~$\B$ and
consider a perturbation $\B\to\B'$ in the class of trigonal
curves. We are interested in the perturbation epimorphism
$\pi_1(\MB\sminus\B)\onto\pi_1(\MB\sminus\B')$, where $\MB$ is a
Milnor ball about~$Q$. We take for~$\MB$ the union of the affine
fibers over a small disk~$\Delta$
about the type~$\tD{_m}$ fiber of~$\B$; then the groups can be
computed using van Kampen's method, \cf\.~\ref{s.presentation}
or~\ref{s.pert.e8}, with only the monodromy within~$\Delta$ taken
into account.

Let~$u$ be the $(m-4)$-valent \cross-vertex of the dessin of~$\B$
representing the singular fiber containing~$Q$. Pick a trivalent
\black-vertex~$v$ adjacent to~$u$ and let $\{\bc_1,\bc_2,\bc_3\}$
be a canonical basis in the fiber~$F_v$ over~$v$, \cf\.
Figure~\ref{fig.basis}, defined by the marking such that $[u,v]$
is the solid edge opposite to~$e_3$ at~$v$. (In other words, the
generators~$\bc_1$ and~$\bc_2$ are brought together when the fiber
approaches~$u$.) Then, the braid monodromy about~$u$ is
$\Gs_1^{m-4}(\Gs_1\Gs_2)^3$,
%; this
and letting $\Gs_1^{m-4}(\Gs_1\Gs_2)^3=\id$
results in the relations
$[\bc_3,\bc_1\bc_2]=1$ and
$\Gs_1^{m-4}(\bc_i)=(\bc_1\bc_2\bc_3)\1\bc_i(\bc_1\bc_2\bc_3)$,
$i=1,2$.
Since the restriction of~$\Gs^{-2}_1$ to the subgroup
$\<\bc_1,\bc_2\>$ is the conjugation by $(\bc_1\bc_2)$, one obtains
$$
\pi_1(\MB\sminus\B)=\<\bc_1,\bc_2\>\rtimes\<\bc_3\>,\quad
\bc_3\1\bc_i\bc_3=\Gs_1^{m-2}(\bc_i),\ i=1,2.
$$

\lemma\label{pert.A}
For the perturbation $\bD_m\to\bA_{m-1}$ \rom(or any further
perturbation thereof\rom), the group $\pi_1(\MB\sminus\B')$ is
abelian.
\endlemma

\proof
The perturbation is realized as follows: the vertex~$u$ is
replaced with a new \cross-vertex~$u'$ of valency~$m$ and the
fragment shown in Figure~\ref{fig.p-dm}\fref(a). (For clarity, we
keep omitting \white-- and \cross-vertices in the drawings; the
new vertex~$u'$ is in the outer region in
Figure~\ref{fig.p-dm}\fref(a). The appearance of the two new
\black-vertices is due to the fact that the perturbation increases
the degree of~$\jB$.) Choosing a new canonical basis
$\{\bb_1,\bb_2,\bb_3\}$ over one of the \black-vertices in the
figure, from the two loops one obtains the relations $\bb_2=\bb_3$
and $\bb_2=\bb_1\bb_2\bb_3\bb_2\1\bb_1\1$. Hence, the group is
abelian.
\endproof

\midinsert
\centerline{\vbox{\halign{\hss#\hss&&\qquad\quad\hss#\hss\cr
\cpic{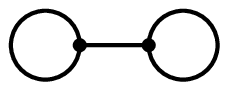}&\cpic{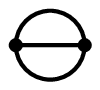}&\cpic{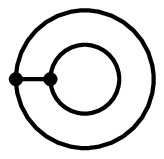}\cr
\noalign{\medskip}
(a) $\bA_{m-1}$&(b) $\bD_2\splus\,\ldots$&(c) $\bD_3\splus\,\ldots$\cr
\crcr}}}
\figure\label{fig.p-dm}
The perturbations of~$\bD_m$
\endfigure
\endinsert

\lemma\label{pert.D}
Consider the perturbation
$$
\tsize
\bD_m\to\bD_p\splus\bigoplus_{i=1}^k\bA_{s_i},\quad
p\ge2,\quad
d:=m-p-\sum_{i=1}^k(s_i+1)\ge0.
$$
If $d=0$, let $s=\gcd(s_i+1)$,
$1\le i\le k$\rom; otherwise, let $s=1$.
Then
$$
\pi_1(\MB\sminus\B')=\TG{s}\rtimes\<\bc_3\>,\quad
\bc_3\1\bc_i\bc_3=\Gs_1^{m-2}(\bc_i),\ i=1,2,
$$
where $\TG{s}=\<\bc_1,\bc_2\,|\,\{\bc_1,\bc_2\}_s=1\>$.
%$\pi_1(\MB\sminus\B')$ is a semidirect product of the group
%$\<\bc_1,\bc_2\,|\,\{\bc_1,\bc_2\}_s=1\>$ and the free
%group~$\<\bc_3\>$, the generator~$\bc_3$ acting on the kernel
%\via\ $\Gs_1^{p-2}$.
\endlemma

\proof
The perturbation is realized as follows: the original
$(m-4)$-valent vertex~$u$ is
replaced with:
\roster
\item\local{pf.1}
$k$ \cross-vertices of valencies
$s_1+1,\ldots,s_k+1$,
\item\local{pf.2}
%$m-p-\sum_{i=1}^k(s_i+1)$
$d$ monovalent \cross-vertices, and
\item\local{pf.3}
either the fragment
shown in Figure~\ref{fig.p-dm}\fref(b) (if $p=2$)
or~\fref(c) (if $p=3$),
or a $(p-4)$-valent \cross-vertex~$u'$ (if $p>4$).
\endroster
(If $p\le3$, all singular fibers are of type~$\tA{}$ and
the perturbation increases the degree
of~$\jB$, introducing two new \black-vertices. If $p>4$, the fiber
over~$u'$ is of type~$\tD{}$. If $p=4$, there also is a
type~$\tD{_4}$ singular fiber of~$\B'$ that does not correspond to
any vertex of the dessin.)

The braid monodromy about the $i$-th vertex in~\loccit{pf.1} is
$\Gs_1^{s_i+1}$, and the
%braid
monodromy about each vertex
in~\loccit{pf.2}, if any, is $\Gs_1$. Thus, the braid relations
resulting from~\loccit{pf.1} and~\loccit{pf.2} simplify to
$\Gs_1^s=\id$, or $\{\bc_1,\bc_2\}_s=1$. The
%braid
monodromy about
the original type~$\tD{_m}$ singular fiber (the `monodromy at
infinity') is
$\Gs_1^{m-4}(\Gs_1\Gs_2)^3$;
as explained at the beginning of this section, it results in the
braid relations
$$
\bc_3\1\bc_i\bc_3=\Gs_1^{m-2}(\bc_i),\quad i=1,2,
\eqtag\label{eq.d-braid}
$$
as stated.

If $p\ge4$, there are no other relations, as the remaining
type~$\tD{_p}$ fiber can be ignored in the presence of the
`relation at infinity'. Otherwise, if $p=2$ or~$3$, the braid
$\Gs_1^{p-4}(\Gs_1\Gs_2)^3$ above is the monodromy along the outer
contour of the insertion shown in Figure \ref{fig.p-dm}. The braid
relations resulting from the two regions separately are
$[\bc_1,\bc_3]=[\bc_2,\bc_3]=1$ (if $p=2$,
Figure~\ref{fig.p-dm}\fref(b)) or
$\{\bc_2,\bc_3\}_4=1$ and
$\bc_2=\bc_3\1\bc_2\1\bc_1\bc_2\bc_3$
(if $p=3$,
Figure~\ref{fig.p-dm}\fref(c));
they are equivalent to~\eqref{eq.d-braid}.
\endproof

\Remark
Formally, one has $\pi_1(\MB\sminus\B)=\TG0\rtimes\Z$
and
$\pi_1(\MB\sminus\B')$ is obtained from
$\pi_1(\MB\sminus\B)$ by adding an extra relation
$\{\bc_1,\bc_2\}_s=1$.
\endRemark

\Remark
If $s$ is odd, %then
$\Gs_1$ is an inner automorphism of~$\TG{s}$ and, in any case, $\Gs_1^2$
is an inner automorphism of~$\TG{s}$. Hence, if $s$ is odd or $m$ is
even, the group in Lemma~\ref{pert.D} splits into \emph{direct} product
$\TG{s}\times\Z$.
\endRemark

\corollary
For a perturbation as in Lemma~\ref{pert.D}, the group
$\pi_1(\MB\sminus\B')$ is abelian if and only if
either $s=1$ or $s=2$ and $m$ is even.
\qed
\endcorollary

\corollary\label{pert.d6}
The only perturbations of a type~$\bD_6$ singular point that have
nonabelian fundamental groups are $\bD_3\splus\bA_2$ and
$\bD_2\splus\bA_3$.
\qed
\endcorollary

\subsection{Perturbations of maximal sextics}
According to the results of~\S\ref{S.groups}, the fundamental
group of a \emph{maximal} sextic~$B$
satisfying~\cstar\
%as in Theorem~\ref{th.main}
is abelian unless the set of singularities of~$B$ is
$\bE_8\splus\bA_4\splus\bA_3\splus2\bA_2$ or
$\bE_8\splus\bD_6\splus\bA_3\splus\bA_2$. In this section, we show
that the only perturbations of these two sextics that have
nonabelian groups are those listed in Theorem~\ref{th.pert}.

\lemma\label{pert.irr}
Let~$B$ be the irreducible plane
sextic with the set of singularities
$\bE_8\splus\bA_4\splus\bA_3\splus2\bA_2$. The only proper
perturbation $B\to B'$ that has nonabelian fundamental group is given
by $\bE_8\to\bA_4\splus\bA_3$\rom; in this case, the perturbation
epimorphism is an isomorphism.
\endlemma

\proof
Any perturbation of the $\bA$-type points of~$B$ can be realized
on the level of the trigonal model~$\B$: in the language of the
skeletons, the \cross-vertex at the center of a region of~$\Sk$
splits into several \cross-vertices of smaller valencies. Assume
that a vertex of valency~$l$ splits into vertices of valencies
$l_1,\ldots,l_k$, so that $l_1+\ldots+l_k=l$. Under an appropriate
choice of the basis $\{\Ga_1',\Ga_2',\Ga_3'\}$, the braid relation
about the original \cross-vertex is $\{\Ga_1',\Ga_2'\}_l=1$.
%Then,
%using~\eqref{eq.equiv}, one can see that after the splitting this
After the splitting, this
relation simplifies to
$\{\Ga_1',\Ga_2'\}_{l'}=1$, where
$l'=\gcd(l_i)$, $1\le i\le k$, see~\eqref{eq.equiv}.

Applying this observation to the
%case under consideration,
set of singularities $\bE_8\splus\bA_4\splus\bA_3\splus2\bA_2$,
one can
see that, if the point perturbed is~$\bA_4$, $\bA_3$, or~$\bA_2$
(so that $l$ above is~$5$, $4$, or~$3$, respectively),
%then
the new
parameter~$l'$ divides~$1$, $2$, or~$1$, respectively.
Hence, instead of $(l,m,n)=(5,4,3)$, see~\ref{s.lmn}, one has
$(l,m,n)=(1,4,3)$, $(5,2,3)$, or~$(5,4,1)$, respectively.
(The two
%type~$\bA_2$ points
cusps
of~$B$ are permuted by the complex
conjugation. Hence, if a cusp is perturbed,
one can assume that it is the one
over~$t$, see Figure~\ref{fig.regions}.)
Using
\GAP~\cite{GAP}, see~\ref{s.others} (for the last case, see
also~\ref{s.stem}), one concludes that all three groups are
abelian.

Assume that perturbed is the
%distinguished
type~$\bE_8$ point~$P$.
Let~$\MB$ be a Milnor ball about~$P$, and let
$\{\bc_1,\bc_2,\bc_3\}$ be the basis for $\pi_1(\MB\sminus B)$
introduced at the beginning of~\ref{s.pert.e8}. As shown
in~\cite{dessin.e7}, the inclusion homomorphism
$\pi_1(\MB\sminus B)\to\pi_1(\Cp2\sminus B)$ is given by
$$
\bc_1\mapsto(\Ga_1\Ga_2)\Ga_3(\Ga_1\Ga_2)\1,\quad
\bc_2\mapsto\Ga_1,\quad
\bc_3\mapsto\Ga_3.
\eqtag\label{eq.incl.e8}
$$
%Note that, in
In particular,
it follows from~\eqref{eq.incl.e8} and~\iref{s.720}{G.6} that
%the inclusion induces
it is
an epimorphism. Hence, if
$\pi_1(\MB\sminus B')$ is abelian, so is $\pi_1(\Cp2\sminus B')$.
In the remaining three cases, one adds to~\eqref{eq.G} the
relations between (the images of) $\bc_1$, $\bc_2$, $\bc_3$ given
by Proposition~\ref{pert.e8} and computes the sizes of the new
groups. The results are $720$ (the group does not change), $6$,
and~$6$ (the group is abelian).
\endproof

\lemma\label{pert.red}
Let~$B$ be the reducible plane
sextic with the set of singularities
$\bE_8\splus\bD_6\splus\bA_3\splus\bA_2$. The only proper
perturbation $B\to B'$ that has nonabelian fundamental group is given
by $\bE_8\to\bD_5\splus\bA_2$\rom; in this case, the perturbation
epimorphism is an isomorphism.
\endlemma

\proof
If the type~$\bA_4$ or type~$\bA_2$ point is perturbed, then, as
explained
in the previous proof, one replaces $(l,m,n)=(4,3,\0)$
in~\eqref{eq.Gi} with $(l,m,n)=(2,3,\0)$ or $(4,1,\0)$,
respectively. Computing the size of the quotient by~$\Ga_2^3$,
one concludes that the group is abelian.

Assume that the type~$\bD_6$ singular point~$Q$ is perturbed and
let~$\MB$ be a Milnor ball about~$Q$. The inclusion homomorphism
$\pi_1(\MB\sminus B)\to\pi_1(\Cp2\sminus B)$ is an epimorphism.
(This is true for any triple point of any trigonal curve, as the
three generators of $\pi_1(\Cp2\sminus B)$, when chosen in a fiber
close to~$Q$, `fit' into~$\MB$.) Thus, if
$\pi_1(\MB\sminus B')$ is abelian, so is $\pi_1(\Cp2\sminus B')$
and, in view of Corollary~\ref{pert.d6}, it remains to consider
the perturbations (back to the conventional notation)
$2\bA_1\splus\bA_3$ and $\bA_3\splus\bA_2$, which result in extra
relations $\{\Ga_2,\Ga_t\}_s=1$, $s=4$ or~$3$, respectively, \cf\.
Lemma~\ref{pert.D}, Figure~\ref{fig.e8-r}\fref(d),
and~\eqref{eq.lmn}. In other words, the values $(l,m,n)=(4,3,\0)$,
see~\ref{s.120},
%should be
are replaced with $(4,3,4)$ or $(4,3,3)$,
respectively. Both groups are abelian.

Finally, the perturbations of the type~$\bE_8$ point~$P$ are
studied similar to the previous proof, using~\eqref{eq.incl.e8}.
In view of~\iref{s.120}{Gi.6}, whenever $\pi_1(\MB\sminus B')$ is
abelian, so is $\pi_1(\Cp2\sminus B')$. For the three nonabelian
perturbations of~$P$, Proposition~\ref{pert.e8},
the sizes of the quotient
$\pi_1(\Cp2\sminus B')/\Ga_1^3$ are $15$, $15$ (the group is
abelian), and $1800$ (the group does not change,
\cf\.~\iref{s.120}{Gi.1}).
\endproof

\subsection{Degenerations of sextics}
In this section, we show that each plane sextic $B$ satisfying~\cstar\
degenerates to a maximal one.

\lemma\label{lem.irr}
Any irreducible plane
sextic~$B$ with a type~$\bE_8$ singular point and
at least two more triple points
%can be obtained by a perturbation from
degenerates to
a maximal sextic with the set of singularities
$\bE_8\splus\bE_6\splus\bD_5$ \rom(No\.~$30$ in
Table~\ref{tab.e8}\rom).
\endlemma

\proof
Consider the triangle Cremona transformation
$\Cp2\dashrightarrow\Cp2$ with the centers at
the type~$\bE_8$ point and
%at
two other triple points
of~$B$. The transform of~$B$ is a cuspidal cubic $C\subset\Cp2$,
and the three exceptional divisors~$E_i$, $i=1,2,3$, are
positioned as follows:
\Dashes
\dash
$E_1$ is tangent to~$C$ at the cusp, and
\dash
none of the intersection points $E_i\cap E_j$, $1\le i<j\le3$,
belongs to~$C$.
\endDashes
It is clear that one can keep~$E_1$ and degenerate the pair
$(E_2,E_3)$ to an `extremal' position, so that, up to reordering,
$E_2$ is inflection tangent to~$C$ and $E_3$ is tangent to~$C$.
(Recall that $C$ has a unique inflection point.)
The new inverse transform of~$C$ has the set of singularities
$\bE_8\splus\bE_6\splus\bD_5$.
\endproof

\lemma\label{lem.red}
Any reducible plane
sextic~$B$ with a type~$\bE_8$ singular point and
at least two more triple points
%can be obtained by a perturbation from
degenerates to
a maximal sextic with the set of singularities
$\bE_8\splus\bD_6\splus\bD_5$
%\rom(No\.~$16'$ in Table~\ref{tab.e8-r}\rom)
or $\bE_8\splus\bE_7\splus\bD_4$
%\rom(No\.~$17'$ in Table~\ref{tab.e8-r}\rom).
\rom(respectively, No\.~$16'$ or~$17'$ in
Table~\ref{tab.e8-r}\rom).
\endlemma

\proof
We proceed as in the previous proof. This time, $B$ splits into an
irreducible quintic~$B_5$ and a line~$B_1$.
The linear component~$B_1$ passes through two double
points $P_2$, $P_3$ of~$B_5$, and we will deform~$B_5$,
keeping~$B_1$ passing through
%(the images of)
these points.
Apply the triangular Cremona transformation
$\Cp2\dashrightarrow\Cp2$
centered at the type~$\bE_8$ point, $P_2$,
and~$P_3$. The transform of~$B_5$ is a cuspidal cubic~$C$, and the
exceptional divisors $E_i$, $i=1,2,3$, are as follows:
\Dashes
\dash
$E_1$ is tangent to~$C$ at the cusp,
\dash
$E_2\cap E_3\in C$, and
\dash
the other two points $E_1\cap E_2$, $E_1\cap E_3$ do \emph{not}
belong to~$C$.
\endDashes
Now, one can keep~$E_1$ and degenerate $(E_2,E_3)$ to one of the
following `extremal' configurations: either $E_2$ is inflection
tangent to~$C$ (necessarily at $E_2\cap E_3$), or $E_2$ is tangent
to~$C$ at $E_2\cap E_3$ and $E_3$ is tangent to~$C$ at another
point. (Recall that $C$ has a unique inflection point~$Q_0$, and
that from each smooth point $Q\ne Q_0$ of~$C$, there is a unique
tangent to~$C$ other than that at~$Q$.) The new inverse transform
of~$C$, with the line $B_1=(P_2,P_3)$ added, has the set of
singularities $\bE_8\splus\bE_7\splus\bD_4$ or
$\bE_8\splus\bD_6\splus\bD_5$, respectively.
\endproof

\proposition\label{->max}
Each plane sextic~$B$
satisfying~\cstar\
degenerates to a maximal sextic
satisfying~\cstar.
\endproposition

\proof
Due to Lemmas~\ref{lem.irr} and~\ref{lem.red}, it suffices to
consider the case when
$B$ has at most one triple point other than~$P$. Then the trigonal
model~$\B$ of~$B$ is not isotrivial and
has at most one triple point; hence, it is
obtained by at most one elementary transformation from a trigonal
curve~$\B'$ with double singular points only. According
to~\cite{degt.kplets}, there is a degeneration $\B'_t$,
$t\in[0,1]$, of $\B'=\B'_1$ to a maximal curve~$\B'_0$; it is
followed by a degeneration (in the class of trigonal models of
sextics satisfying~\cstar)
$\B_t$ of $\B=\B_1$ to a maximal
trigonal model~$\B_0$ and, hence, by a degeneration $B_t$ of
$B=B_1$ to a maximal sextic~$B_0$.
\endproof

\Remark
If $B$ has three triple points then, in the proof of
Proposition~\ref{->max}, the trigonal model~$\B$ is obtained from
$\B'\subset\Sigma_1$ by two elementary transformation and, during
the degeneration $\B'_t$, the two fibers contracted could merge to
a single fiber, resulting in a non-simple singular point of~$\B$.
To exclude this possibility, we treat the case of three triple
points separately in Lemmas~\ref{lem.irr} and~\ref{lem.red}.
\endRemark

\subsection{Proof of Theorem~\ref{th.main}}\label{proof.main}
According to Proposition~\ref{->max}, any plane sextic
satisfying~\cstar\ degenerates to a maximal one; hence, due to the
calculation in~\S\ref{S.groups}, nonabelian can only be the groups
of the perturbations of the sextics with the sets of singularities
$\bE_8\splus\bA_4\splus\bA_3\splus2\bA_2$ or
$\bE_8\splus\bD_6\splus\bA_3\splus\bA_2$. (Recall that a
perturbation $B\to B'$ induces an epimorphism
$\pi_1(\Cp2\sminus B)\onto\pi_1(\Cp2\sminus B')$,
see~\cite{Zariski.group}.) According to
Lemmas~\ref{pert.irr} and~\ref{pert.red}, there are only two
proper perturbations with nonabelian fundamental groups; none of
them has a type~$\bE_8$ singular point.
\qed

\subsection{Proof of Theorem~\ref{th.pert}}\label{proof.pert}
The statement follows immediately from the list of sextics with
nonabelian groups (Theorem~\ref{th.main}) and the description of
their perturbations
(Lemmas~\ref{pert.irr} and~\ref{pert.red}).
\qed

\widestnumber\key{EO1}
\refstyle{C}

\enddocument